\input amstex
\loadeufm
\loadeufb
\loadeurm
\loadeurb
\loadeusm
\loadeusb
\documentstyle{amsppt}
\newcount\nagyitas
\nagyitas=1

\define\get#1{\input $HOME/tex/amstex/#1}
\define\Title{}

\define\shorttitle{\Title}
\newcount\coauth
\coauth=0
\newcount\abst
\abst=0

\newcount\cim
\cim=1
\def\NoTitle{\cim=0}
\newcount\ded
\ded=0

\newcount\scl
\scl=0

\def\today{\ifcase \the\month\or January \or February \or March \or
April \or May \or June \or July \or August \or September \or October
\or November \or December \fi \the\day, \the\year}

\def\mypagesize{
\pageheight{18.5cm}
\pagewidth{14cm}
}

\newcount\specialnum
\specialnum=0

\newcount\whospage
\whospage=0

\scl=0
\NoTitle

\redefine\Title
{\uppercase {
Logarithmic Kodaira-Akizuki-Nakano vanishing for singular varieties
}} 
\redefine\shorttitle{\uppercase {
Logarithmic Kodaira-Akizuki-Nakano vanishing
}}

%
\pageheight{20cm}
\vcorrection{-.25cm}
\pagewidth{14.5cm}
\hcorrection{-.25cm}

\topmatter
\phantom{m}\vskip-.5cm
\title 
Logarithmic Kodaira-Akizuki-Nakano vanishing and  
Arakelov-Parshin boundedness for singular varieties
\endtitle 
\author 
S\'andor J. Kov\'acs 
\endauthor
\address
Department of Mathematics, 
University of Chicago, Chicago, IL 60637
\endaddress
\email {\tt 
skovacs\@member.ams.org} 
\endemail 
\date \today \enddate
\thanks
Supported in part by an AMS Centennial Fellowship and by NSF Grant
DMS-9818357.
\endthanks
\endtopmatter
\rightheadtext{\shorttitle}

\ifnum\whospage=0
\mypagesize
\fi

\ifnum \nagyitas > 0
\magnification=\magstep1
\fi


\define\noin{\noindent}
 2
 4
 2

 1

\newcount\thm
\newcount\sect
\newcount\subsect
\newcount\oldthm
\newcount\oldsect
\newcount\whatis
\newcount\nothmscount
\thm=0 \sect=0 \oldsect=0 \subsect=0\whatis=0 \nothmscount=0
\define\thmnumber{\the\sect.\the\thm}
\define\nosections{\redefine\thmnumber{\the\thm}}
\define\section#1{\bigskip \ifnum\whatis=1\sect=\oldsect 
\redefine\thmnumber{\the\sect.\the\thm}\fi\whatis=0
\thm=0 \advance\sect by1 \heading \S 
\the\sect.\  #1 \endheading 
} 

\define\subsection#1{\advance \thm by1\advance \oldthm by1
\ifnum \whatis=0 \vskip-6pt\oldthm=\thm\oldsect=\sect\sect=\thm 
        \else \medskip \sect=\oldthm \fi
\whatis=1
\thm=0  
\redefine\thmnumber{\the\oldsect.\the\sect.\the\thm}
	\ifnum\nothmscount=1
	\redefine\thmnumberfortags{\the\oldsect.\the\sect}  
	\else
	\redefine\thmnumberfortags{\the\oldsect.\the\sect.\the\thm}  
	\fi
\subheading{\S\S 
\the\oldsect.\the\sect\  #1}
\noindent}
\define\endsubsection{
\ifnum\whatis=1\sect=\oldsect
\advance \thm by1\niltag
\redefine\thmnumber{\the\sect.\the\thm}
\redefine\thmnumberfortags{\the\sect.\the\thm}
\fi\whatis=0\nothmscount=0
}
\settabs \+ \hskip .5in & \cr
\newcount\rostit
\define\myroster{\begingroup \parindent=48pt}
\define\endmyroster{\par\relax\endgroup}
\define\myitem{
\item{$(\thmnumber.\the\mytagcounter)$}
{\global\advance\mytagcounter by1}
} 
\def\rostlabel#1{{\advance\mytagcounter
by-1\taglabel#1\advance\mytagcounter by1}}
\define\myfootnote#1{\parindent=12pt\footnote{#1}\parindent=48pt} 
\def\label#1{ 
\edef\temp{\thmnumber}
\global\let #1 = \temp }
\define\quote#1{$(#1)$}
\def\itemlabel#1{
\edef\tempitem{\thmnumber.\the\rostit}
\global\let #1 = \tempitem}
\TagsOnRight
\CenteredTagsOnSplits
\newcount\mytagcounter
\mytagcounter=1
\define\thmnumberfortags{\thmnumber}
\define\niltag{\global\mytagcounter=1}
\define\mytagcount{\thmnumberfortags.\the\mytagcounter}
\def\taglabel#1{
\edef\ttemp{\thmnumberfortags.\the\mytagcounter}\global\let#1=\ttemp}
\define\mytag#1{
\taglabel{#1}\tag{#1}
{\global\advance\mytagcounter by1}}

\newskip\skegy
\newskip\skket
\newskip\skhar
\newskip\sknegy
\newskip\skot
\define\mytab#1 #2 #3 #4 {\begingroup
\skegy=#1cm
\skket=#2cm
\skhar=#3cm
\sknegy=#4cm
\skot=\pagewidth
\cleartabs 
\tabalign \hskip\skegy & \hskip\skket &
\hskip\skhar & \hskip\sknegy &\cr
\vskip -\baselineskip
\advance\skot by -\skegy
\advance\skot by -\skket
\advance\skot by -\skhar
\advance\skot by -\sknegy
\advance\skket by -.2cm
\advance\skhar by -.2cm
\advance\sknegy by -.2cm
}

\define\tabitem#1#2#3#4{\tabalign 
&\hbox{\hsize=\skket\vtop{\noin #1}}
&\hbox{\hsize=\skhar\vtop{\noin #2}} 
&\hbox{\hsize=\sknegy\vtop{\noin #3}} 
&\hbox{\hsize=\skot\vtop{\noin #4}}
\cr
}
\ifnum\specialnum=0
\define\num{\global\advance\thm by1 \niltag{\thmnumber\ }}
\define\remnum{\global\advance\thm by0 \niltag{\thmnumber.1\ }}
\else
\define\numero{\global\advance\thm by1 \niltag{\thmnumber\ }}
\define\remnum{\global\advance\thm by0 \niltag{\thmnumber.1\ }}
\fi
\define\subnum{\noin
{\mytagcount}
{\global\advance\mytagcounter by1}}
\define\newnum{\global\advance\thm by1 \niltag\demo{\bf \thmnumber}}
\define\endnewnum{\enddemo}
\define\defini{\global\advance\thm by1 \niltag\demo{\smc \thmnumber\
Definition}} 

\define\proof{\demo{\bf Proof}}
\define\endproof{$\square$\enddemo}


\define\kan{Kodaira-Akizuki-Nakano vanishing theorem}
\define\moc{morphism of complexes}

\define\dt{distinguished triangle}
\define\te{there exist}
\define\st{such that}

\define\cf#1{cf.\ \cite{#1}}
\define\shes{short exact sequence}


\define\m{\eusm M}

\redefine\l{\eusm L}

\define\C{\Bbb {C}}
\define\N{\Bbb {N}}

\define\Q{\Bbb {Q}}



\def\tensor{\otimes}

\def\rpforward#1{R{#1}_*}
 
\def\functor{\Gamma}

\define\opname{\operatorname}
\define\id{\operatorname{id}}

\redefine\hom{\operatorname{Hom}}
\ifnum\specialnum=0

\else

\fi

\define\w#1{\bigwedge^#1}
\define\fundgp#1{\pi_1(#1,\star)}
\define\pic #1{\operatorname{Pic}(#1)}
\define\pico #1{\operatorname{Pic}^{\circ}(#1)}
\define\aut#1{\operatorname{Aut}(#1)}
\def\sym{\operatorname{Sym}}
\def\rank{\operatorname{rk}}

\def\im{\operatorname{im}}
\def\ob{\operatorname{Ob}}

\def\supp{\operatorname{supp}}
\redefine\ne#1{\overline{NE}(#1)}

\def\Spec{\operatorname{Spec}}
\def\spec#1.#2.{{\bold S\bold p\bold e\bold c}_{#1}#2}
\def\ring#1.{\Cal O_{#1}}
\def\map#1.#2.{#1 \to #2}
\def\longmap#1.#2.{#1 \longrightarrow #2}
\def\pe#1.{\Bbb P(#1)}
\def\pr#1.{\Bbb P^{#1}}
\def\proj#1.#2.{{\bold P\bold r\bold o\bold j}_{#1}\sum #2}

\define\coh#1.#2.#3.{H^{#1}(#2,#3)}
\define\dimcoh#1.#2.#3.{h^{#1}(#2,#3)}
\define\hypcoh#1.#2.#3.{\Bbb H_{\vphantom{l}}^{#1}(#2,#3)}
\define\loccoh#1.#2.#3.#4.{H^{#1}_{#2}(#3,#4)}
\define\dimloccoh#1.#2.#3.#4.{h^{#1}_{#2}(#3,#4)}
\define\lochypcoh#1.#2.#3.#4.{\Bbb H^{#1}_{#2}(#3,#4)}
\define\ses#1.#2.#3.{0  \longrightarrow  #1   \longrightarrow 
 #2  \longrightarrow   #3 \longrightarrow  0}
\define\sesshort#1.#2.#3.{0  @>>>  #1  @>>>  #2  @>>>  #3 @>>> 0}
\define\dist#1.#2.#3.{  #1   \longrightarrow 
 #2  \longrightarrow   #3 @>+1>> }  
\define\CDdist#1.#2.#3.{  #1   @>>>  #2  @>>>   #3 @>+1>> }  
\define\shortses#1.#2.#3.{0  \rightarrow  #1   \rightarrow 
 #2  \rightarrow   #3 \rightarrow  0}
\define\shortdist#1.#2.#3.{  #1   \rightarrow 
 #2  \rightarrow   #3 @>+1>> }  
\def\ddist#1.#2.#3.#4.#5.#6.{\CD
#1 @>>> #2 @>>> #3 @>+1>> \\
@VVV @VVV @VVV \\
#4 @>>> #5 @>>> #6 @>+1>> 
\endCD}
\def\ddistun#1.#2.#3.#4.#5.#6.{\CD
#1 @>>> #2 @>>> #3 @>+1>> \\
@. @VVV @VVV  \\
#4 @>>> #5 @>>> #6 @>+1>> 
\endCD}
\define\factor#1.#2.{\left. \raise 2pt\hbox{$#1$} \right/
\hskip -2pt\raise -2pt\hbox{$#2$}}
\def\iff#1#2#3{
\hfil\hbox{\hsize =#1
\vtop{\noin #2}
\hskip.5cm 
\lower.5\baselineskip\hbox{$\Leftrightarrow$}\hskip.5cm
\vtop{\noin #3}}\hfil\medskip}
\define\myoplus#1.#2.{\underset #1 \to {\overset #2 \to \oplus}}

\define\oldresto#1{\hskip -.5pt\hbox{$\mid_{#1}$}}
\define\resto#1{\hbox{\hbox{$\big\vert_{#1}$}}}
\define\qis{\,{\simeq}_{qis}\,}
\define\wedp{\wedge_{p}}
\define\ww#1.#2.{\curlywedge
_{#1}^{#2}}
\define\wws#1.#2.{\scriptstyle\ww#1.#2.}

\define\fh#1.#2.{F^{#1} R^{#2}}
\define\epq#1.#2.#3.{E^{#1,#2}_{#3}}
\define\dd#1.#2.#3.{d^{\, #1,#2}_{#3}}
\define\kk#1.#2.#3.{K^{\, #1,#2}_{#3}}
\define\ki#1.#2.#3.{I^{\, #1,#2}_{#3}}
\define\ff#1.#2.{\frak F_{#1}^{#2}}
\define\ffb#1.#2.{\boxed{\frak F_{#1}^{#2}}}
\define\ffbx#1.#2.#3.{\boxed{\frak F_{#1}^{#2}(#3)}}
\define\ffs#1.#2.{\scriptstyle\frak F_{#1}^{#2}}
\define\ffsb#1.#2.{\boxed{\ffs#1.#2.}}
\define\ffsbx#1.#2.#3.{\boxed{\scriptstyle\frak F_{#1}^{#2}(#3)}}
\define\fa#1.#2.{\frak f_{#1}^{#2}}
\define\fab#1.#2.{\boxed{\frak f_{#1}^{#2}}}
\define\ffa#1.#2.{\frak F\frak i\frak l\frak t_{#1}^{#2}}
\define\ffab#1.#2.{\boxed{\frak F\frak i\frak l\frak t_{#1}^{#2}}}
\define\rfi#1.#2.{R^{#1}\functor(#2\,)}
\define\hypf#1{\Bbb F^{#1}_{\phantom{\cdot}}}
\define\hypg#1{\Bbb G^{#1}_{\phantom{\cdot}}}
\define\hype#1{\Bbb E^{#1}_{\phantom{\cdot}}}
\define\hypgre#1{\Bbb Gr^{#1}_{\Bbb E}}
\define\hypgr#1{\Bbb Gr^{#1}_{\Bbb F}}
\define\hypgrr#1#2{\Bbb Gr^{#1}_{#2}}
\define\al#1.#2.#3.{\alpha_{#1,#2}^{#3}}
\define\be#1.#2.{\beta_{#1,#2}}
\define\ga#1.#2.{\gamma_{#1,#2}}
\define\yy#1.#2.#3.#4.{y_{#1,#2}^{#3}(#4)}
\define\zz#1.#2.#3.#4.{z_{#1,#2}^{#3}(#4)}
\def\\{\char`\\}

\pageheight{20cm}
\vcorrection{-.25cm}
\pagewidth{14.45cm}
\hcorrection{-.25cm}

\input diagrams.tex

\diagramstyle[amstex]
\define\Ox#1.#2.{\underline{\Omega}_X^{#1}(\log{#2})}
\define\Oxx#1.#2.{K_{#1}}
\define\Oxp{\underline{\Omega}_{X}^{p}(\log{D})}
\define\Oxcp{\underline{\Omega}_{X/C}^{p}(\log{D})}
\define\Oxc#1.{\underline{\Omega}_{X/C}^{#1}(\log{D})}
\define\Oc#1.{{\omega}_C^{}({#1})}
\define\Oy#1.#2.{\underline{\Omega}_Y^{#1}(\log{#2})}
\define\Oyc#1.#2.{\underline{\Omega}_{Y/C}^{#1}(\log{#2})}
\define\Ol#1.#2.{\underline{\Omega}_L^{#1}(\log{#2})}

\define\letgoodhyp{$\varepsilon_\cdot: (X_\cdot, D_\cdot)\to (X, D)$
 be a good hyperresolution}
\define\tr{\opname{Tr_L}}
\
\define\ygen{{Y_{\text{gen}}}}
\define\gyc{h_*\omega_{Y/C}^m}
\define\fxc{f_*\omega_{X/C}^m}

\newarrow{Sub}{C}{}{}{}{}
\document 


Vanishing theorems have played a central role in algebraic geometry,
for the last couple of decades, especially in classification theory.
\cite{Koll\'ar87} gives an introduction to the basic use of
vanishing theorems as well as a survey of results and applications
available at that time. For more recent results one should consult
\cite{Ein97, Koll\'ar97, Kov\'acs00c, Smith97}. Because of the
availability of such surveys here we will only recall statements that
are important to this article.
In any case one must start with the fundamental vanishing theorem of
Kodaira:

\proclaim{\num Theorem}\cite{Kodaira53} 
Let $X$ be a smooth complex projective variety and $\l$ an ample line
bundle on $X$. Then $$\coh i.X.\omega_X\tensor\l.=0\qquad\text{for
}i>0.$$
\endproclaim\label\kodaira

This has been generalized in several ways, but as noted above we will
restrict to a select few. 
Akizuki and Nakano extended Kodaira's vanishing to include other
exterior powers of the sheaf of differential forms:

\proclaim{\num Theorem}\cite{Akizuki-Nakano54} 
Let $X$ be a smooth complex projective variety and $\l$ an ample line
bundle on $X$. Then $$\coh q.X.\Omega^p_X\tensor\l.=0\qquad\text{for
}p+q>n.$$
\endproclaim\label\kanlabel

The original statement of Kodaira was generalized in a different
direction to allow semi-ample and big line bundles in place of ample
ones by Grauert and Riemenschneider:

\proclaim{\num Theorem}\cite{Grauert-Riemenschneider70} 
Let $X$ be a smooth complex projective variety and $\l$ a semi-ample
and big line bundle on $X$. Then $$\coh
i.X.\omega_X\tensor\l.=0\qquad\text{for }i>0.$$ 
\endproclaim\label\gr

 \demo{\remnum Remark} Later ``semi-ample'' was replaced by ``nef'' in
the statement by Kawamata and Viehweg \cite{Kawamata82,
Viehweg82}. \enddemo

\newnum
Ramanujam gave a simplified proof for \quote\kanlabel\ and showed that it
does not hold if one only requires $\l$ to be semi-ample and big
\cite{Ramanujam72}.
\endnewnum

\newnum
On the other hand Navarro-Aznar {\it et al.} proved a version of the
\kan\ for singular varieties that actually implies \quote\kodaira,
\quote\kanlabel\ and \quote\gr\ \cf{Navarro-Aznar88} in \cite{GNPP88}.
\endnewnum\label\navarro

Another kind of generalization was proved by Esnault and Viehweg:

\proclaim{\num Theorem}\cite{Esnault-Viehweg92, 6.4} 
Let $X$ be a smooth complex projective variety and $\l$ an ample line
bundle on $X$. Further let $D$ be a normal crossing divisor on $X$.
Then $$\coh q.X.\Omega^p_X(\log D)\tensor\l.=0\qquad\text{for
}p+q>n.$$
\endproclaim\label\ev

One of the main goals of the present article is to prove a common
generalization of \quote\kodaira, \quote\kanlabel, \quote\gr,
\quote\navarro\ and \quote\ev. For the exact statement we will first
need to introduce some definitions, so it will only be given in (4.3).
For now let us state it in a very informal form:

\proclaim{\num Main Theorem}
The logarithmic Kodaira-Akizuki-Nakano vanishing theorem of Esnault
and Viehweg admits a good generalization to singular varieties.
\endproclaim

One could ask why we need such a generalization. I believe it is an
interesting result on its own. This seems to be supported by the
enthusiasm that greeted \quote\navarro\ (\cf{Steenbrink85}).  On the
other hand it could be viewed as a ``poor man's version'' of the
logarithmic Kodaira-Akizuki-Nakano vanishing theorem for semi-ample
and big line bundles on smooth varieties. Considering that the obvious
generalization fails, this might be the best one can hope for. As an
easy corollary, we also obtain a local version of this vanishing
theorem.

Nevertheless, my original motivation was an actual application. This
theorem is the corner stone of the proof of an Arakelov-Parshin type
boundedness result.  
That result is presented as an application of the Main Theorem,
although it would merit to be called a ``Main Theorem'' itself.
The first interesting consequence of the Main Theorem is a vanishing
theorem for smooth varieties, (6.4). Note that in order to prove it
one has to go through the singular version. (6.4) is a generalization
of \cite{Kov\'acs97, 1.1} and similar in nature to
\cite{Bedulev-Viehweg00, 2.2}.

Next let us take a brief tour of some related problems, and let us
start by recalling that a family of projective curves is called {\it
isotrivial\/} if all but finitely many members of the family are
isomorphic to a fixed curve.

\newnum
Fix $C$, a smooth projective curve of genus $g$ over an algebraically
closed field of characteristic $0$, $\Delta\subset C$ a finite subset
and $q>1$ positive integer. Let $\delta=\#\Delta$.

Shafarevich conjectured at the 1962 ICM in Stockholm that the set,
$\frak S$, of non-isotrivial families of smooth projective curves of
genus $q$ over $C\setminus \Delta$ is finite.
This was confirmed by \cite{Parshin68} for the case $\Delta=\emptyset$
and by \cite{Arakelov71} in general. Their method was to divide the
problem into two parts:
\roster
\item
``Boundedness'': There are only finitely many deformation types of
families in $\frak S$.
\item
``Rigidity'': There are no non-trivial deformations within $\frak S$.
\endroster

The basic question now is whether Shafarevich's conjecture holds in
higher dimensions.
\endnewnum

\newnum 
It is actually more convenient to work with a compactification of the
family, understanding that later we are free to alter it over
$\Delta$. Hence from now on by a family we will mean a non-isotrivial
family over a compact curve.

Considering families over a compact base curve leads to a related
problem. Namely one could ask what can be said about the singular
fibers of the family. On the simplest level, how many are there? In
fact Szpiro did ask this: Is there a lower bound on the number of
singular fibers if $C\simeq \Bbb P^1$?

 \cite{Beauville81} gave the following answer: there are always at
 least 3 singular fibers and there are families with exactly 3. In
 fact Beauville's proof also shows that there is at least 1 singular
 fiber if the base curve is elliptic. In short $2g-2+\delta>0$.

Note that Kodaira surfaces show that there are families over high
 genus curves without any singular fibers.

More recently \cite{Catanese-Schneider95} asked if the same is true
with higher dimensional fibers, and the conjecture of
\cite{Shokurov97} translates to the same: Is it true that for a
familily of varieties of general type, $2g-2+\delta>0$, or
equivalently: Is $\delta\geq 3$ if $g=0$ and $\delta\geq 1$ if $g=1$?

It is interesting to note the wide range of applications this question
relates to: \cite{Catanese-Schneider95} wanted to use this to obtain
good estimates for the size of the automorphism group of a variety of
general type, while \cite{Shokurov97} needed it for proving
quasi-projectivity of certain moduli spaces.
\endnewnum

\newnum
The basic phylosophy of proving boundedness is first proving that for
a family, $f:X\to C$, \te s some $m\gg 0$, \st\
$\deg(f_*\omega^m_{X/C})$ is bounded in terms of the fixed data and
$m$.  The next step then is to use this bound and an appropriate
moduli space (if it exists) to prove boundedness.

The first step will be called ``weak boundedness''. In practice one
proves weak boundedness and if the appropriate moduli space exists,
then boundedness is almost automatic. The step from weak boundedness
to boundedness is an independent question: it is basically the problem
of existence of moduli spaces.
\endnewnum

\newnum
The following is a select list of results related to these questions.

 \cite{Faltings83} studied the Shafarevich problem for families of
abelian varieties and proved that boundedness holds, while rigidity
fails in general.

 \cite{Migliorini95} showed that for families of minimal surfaces
 $\delta\geq 1$ if $g\leq 1$, and \cite{Kov\'acs96} showed the same for
 families of minimal varieties of arbitrary
 dimension. \cite{Kov\'acs97} settled the question for families of
 minimal surfaces and \cite{Kov\'acs00a} for families of canonically
 polarized varieties: In both cases $2g-2+\delta>0$.

 \cite{Oguiso-Viehweg00} proved the same for families of elliptic
 surfaces. Their work completes the case of families of minimal
 varieties of non-negative Kodaira dimension.

 \cite{Bedulev-Viehweg00} proved that boundedness holds for families
 of surfaces of general type and that weak boundedness (and in some
 cases boundedness) holds for families of canonically polarized
 varieties. As a byproduct of their proof they also obtained that
 $2g-2+\delta>0$ in these cases.
\endnewnum

In this article we obtain results regarding both questions. In fact a
simple observation yields that these questions are in fact strongly
related.

\proclaim{\num Theorem}
Weak boundedness implies that $2g-2+\delta>0$. 
\endproclaim\label\weakbound

The main result of the second part of the article is the following. It
is again in an informal form. For more precise statements see $(7.7),
\ (7.9),\ (7.12)$.

\proclaim{\num Theorem}
Fix $C$, $\Delta\subset C$.
Then weak boundedness holds for families of canonically polarized
varieties with rational Gorenstein singularities (or equivalently
singularities that appear on the canonical models of varieties of
general type) over $C\setminus \Delta$ with fixed Hilbert polynomial
admitting a simultaneous resolution of singularities.  In particular
$2g-2+\delta>0$ for these families. Furthermore in some cases
boundedness holds as well.
\endproclaim

As a corollary one obtains weak boundedness for non birationally
isotrivial families of minimal varieties of general type.

A few days before the completion of this article I learnt that
 \cite{Viehweg-Zuo00} proves that $2g-2+\delta>0$ holds for non
 birationally isotrivial smooth families of minimal varieties. As a
 byproduct of their proof they also obtain weak boundedness for these
 families.

\demo{\smc Definitions and Notation} \nopagebreak
Throughout the article the groundfield will always be $\C$, the field
of complex numbers. A {\it complex scheme \/} (resp.\ {\it complex
variety\/}) will mean a separated scheme (resp.\ variety) of finite
type over $\C$.

A divisor $D$ is called $\Bbb Q$-Cartier if $mD$ is Cartier for some
$m>0$. $D$ is called {\it big\/} if $X$ is proper and $|mD|$ gives a
birational map for some $m > 0$ and it is called {\it nef\/} if
$D.C\geq 0$ for every proper curve $C\subset X$. In particular ample
implies nef and big. If $A$ and $B$ are effective divisors, then
$A\cup B$ will denote $\supp(A+B)$.

Let $U$ be an open subset of $X$.  A line bundle $\l$ on $X$ is called
{\it semi-ample with respect to $U$\/} if some positive power of $\l$
is generated by global sections over $U$. It is called {\it
semi-ample\/} if it is semi-ample with respect to $X$. Similarly $\l$
is called {\it ample with respect to $U$\/} if the global sections of
some positive power of $\l$ define a rational map, that is an
embedding on $U$.

A locally free sheaf $\Cal E$ on a scheme $X$ is called {\it
semi-positive\/} (resp.\ {\it ample\/}) if for every smooth complete
curve $C$ and every map $\gamma : C \to X$, any quotient bundle of
$\gamma^{*}\Cal E$ has non-negative (resp.\ {\it positive\/}) degree.


Let $f:X \to S$ be a morphism of schemes. Then $X_s$ denotes the fibre
of $f$ over the point $s\in S$ and $f_s$ denotes the restriction of
$f$ to $X_s$. More generally, for a morphism $\sigma:\map Z.S.$, let
$f_Z : \map X_Z=X\times_S Z.Z.$. If $f$ is composed with another
morphism $g:S \to T$, then for a $t\in T$, $X_t$ denotes the fibre of
$g\circ f$ over the point $t$, i.e., $X_t = X_{S_t}$. $f_Z$ and $X_Z$
may also be denoted by $f_\sigma$ and $X_\sigma$ respectively.

A singularity is called {\it Gorenstein\/} 
if its local ring is a Gorenstein 
ring. A variety is {\it Gorenstein\/} 
if it admits only Gorenstein 
singularities.
Let $X$ be a normal variety and $f :Y \rightarrow X$ a
resolution of singularities. $X$ is said to have {\it rational
singularities\/} if $R^if_*\ring Y.=0$ for all $i>0$.

Let $X$ be a complex scheme of dimension n. $D_{filt}(X)$ denotes the
derived category of filtered complexes of $\Cal O_{X}$-modules with
differentials of order $\leq 1$ and $D_{filt, coh}(X)$ the subcategory
of $D_{filt}(X)$ of complexes $K^{\cdot}$, such that for all $i$, the
cohomology sheaves of $Gr^{i}_{filt}K^{\cdot}$ are coherent (cf.\
\cite{DuBois81}, \cite{GNPP88}).  $D(X)$ and $D_{coh}(X)$ denotes the
derived categories with the same definition except that the complexes
are assumed to have the trivial filtration. The superscripts $+, -, b$
carry the usual meaning (bounded below, bounded above, bounded).
$C(X)$ is the category of complexes of $\Cal O_{X}$-modules with
differentials of order $\leq 1$ and for
$u\in\operatorname{Mor}(C(X))$, $M(u)\in\ob(C(X))$ denotes the mapping
cone of $u$ (cf.\
\cite{Hartshorne66}). Isomorphism in these categories is 
denoted by $\qis$. If $K^{\cdot}$ is a complex in any of the above
defined categories, then $h^i(K^{\cdot})$ denotes the $i$-th cohomolgy
sheaf of $K^{\cdot}$. In particular every sheaf is naturally a complex
with $h^i=0$ for $i\neq 0$.

The right derived functor of an additive functor $F$, if exists, is
denoted by $RF$ and $R^iF$ stands for $h^i\circ RF$. In particular
$\Bbb H^i$ denotes $R^i\Gamma$ and $\Bbb H^i_Z$ denotes $R^i\Gamma_Z$
and $\Cal H^i_Z$ denotes $R^i\Cal H_Z$ where $\Gamma$ is the functor
of global sections and $\Gamma_Z$ is the functor of global sections
with support in the closed subset $Z$ and $\Cal H_Z$ is the functor of
local sections with support in the closed subset $Z$. Note that
according to this terminology if $\phi:Y\to X$ is a morphism and $\Cal
F$ is a coherent sheaf on $Y$, then $R\phi_*\Cal F$ is the complex
whose cohomology sheaves give the usual higher direct images of $\Cal
F$. The derived functor of $\otimes$ is denoted by $\otimes^L$.


The dimension of the empty set is $-\infty$
\enddemo

\subheading{Acknowledgement}
The lion's share of the work related to this article was done while I
enjoyed the hospitality of RIMS at Kyoto University. I would like to
express my warm thanks to members, visitors and staff of RIMS for
making my stay fruitful and memorable. I would like to thank
especially Shigefumi Mori for making my visit possible and also for
the interest he took in my work. I would also like to thank him,
Stefan Helmke, Stefan Kebekus, Yoichi Miyaoka and Noboru Nakayama for
helpful remarks and discussions. 

I am also indepted to Osamu Fujino for saving me from some simple but
crucial errors.


Commutative diagrams were drawn with the help of Paul Taylor's \TeX\
macro package.

\section{De Rham-Du Bois complexes}
In order to state our generalized version of the \kan, we need Du
Bois's generalized De Rham complex.

The original construction of Du Bois's complex, $\Ox\cdot.D.$, is
based on simplicial resolutions. The reader interested in the details
is referred to the original article \cite{DuBois81}. Note also that a
simplified construction was later obtained by \cite{GNPP88} via the
general theory of cubic resolutions.  An easily accessible
introduction can be found in \cite{Steenbrink85}.

The word ``hyperresolution'' will refer to either simplicial or cubic
resolution. Formally the construction of $\Ox\cdot.D.$ is the same
regardless which resolution is used and no specific aspects of either
resolution will be used.

\demo{\bf \num Definition}
Let $X$ be a complex scheme and $D$ a closed subscheme whose
complement is dense in $X$. Then $(X_\cdot, D_\cdot)\to (X, D)$ is a
{\it good hyperresolution} if $X_\cdot\to X$ is a hyperresolution, and
if $U_\cdot=X_\cdot\times_X (X\setminus D)$ and
$D_\cdot=X_\cdot\setminus U_\cdot$, then $D_i$ is a divisor with
normal crossings on $X_i$ for all $i$.
\enddemo\label\hyper

\proclaim{\num Theorem}
\cite{DuBois81, 6.3, 6.5}
Let $X$ be a proper complex scheme of finite type and $D$ a closed
subscheme whose complement is dense in $X$. Then there exists a unique
$\Ox\cdot.D. \in \ob(D_{filt}(X))$ with the following properties,
using the notation: $$\Ox p.D.:=Gr^{p}_{filt}\, \Ox\cdot.D.[p].$$

\myroster
\myitem
Let $j:X\setminus D\to X$ be the inclusion map. Then
$$\Ox\cdot.D. \qis Rj_*\C_{X\setminus D}.$$
\rostlabel\cc

\myitem 
It is functorial, i.e., if $\phi :Y\rightarrow X$ is a morphism of
proper complex schemes of finite type, then there exists a natural map
$\phi^{*}$ of filtered complexes
$$
\phi^{*}:\Ox\cdot.D. \rightarrow R\phi_{*}\Oy{\cdot}.\phi^*D. 
$$
Furthermore, $\, \Ox\cdot.D. \in \ob(D^{b}_{filt, coh}(X))$ and if
$\phi$ is proper, then $\phi^{*}$ is a morphism in $D^{b}_{filt,
coh}(X)$.
\rostlabel\pushmap

\myitem
Let $U \subseteq X$ be an open subscheme of $X$. Then
$$\Ox\cdot.D.\resto U \qis\underline{\Omega}^{\cdot}_U(\log D\resto
U).$$

\myitem
There exists a spectral sequence degenerating at $E_1$ and abutting to
the singular cohomology of $X\setminus D$:
$$
E_1^{pq}={\Bbb H}^q(X, \Ox p.D.) \Rightarrow
H^{p+q}(X\setminus D, \Bbb C).
$$
\rostlabel\hodge

\myitem
If $\varepsilon_\cdot: (X_\cdot, D_\cdot)\to (X, D)$ is a good 
hyperresolution, then $$  
\Ox\cdot.D.\qis R{\varepsilon_\cdot}_*\Omega^\cdot_{X_\cdot}(\log
D_\cdot).$$ In particular $h^i(\Ox p.D.)=0$ for $i<0$.
\rostlabel\dbconstr

\myitem
There exists a natural map, $\ring X.\to \Ox 0.D.$,
compatible with \quote\pushmap.
\rostlabel\zeromap

\myitem
If $X$ is smooth and $D$ is a normal crossing divisor, then
$$\Ox\cdot.D.\qis\Omega^\cdot_X(\log D).$$
In particular 
$$\Ox p.D.\qis\Omega^p_X(\log D).$$
\rostlabel\dbsmooth

\myitem
If $\phi:Y\to X$ is a resolution of singularities, then
$$ \Ox\dim X.D.\qis R\phi_*\omega_Y(\phi^* D).$$
\rostlabel\topomega

\endmyroster
\endproclaim\label\dub

\section{A short exact sequence}
\noin
The following notation and assumptions will be used throughout this
and the next section.

\newnum Let $X$ be a projective variety and 
$D\subset X$ an effective
divisor on $X$ and
$\varepsilon_\cdot:(X_\cdot, D_\cdot)\to (X, D)$ a good
hyperresolution.
Let $\m$ be a semi-ample line bundle on $X$.
Assume that $\m$ is ample with respect to $X\setminus D$. Let
$\l=\m^N$ for some $N\gg 0$, $\sigma\in \coh 0.X.\l.$ a general
section, and $L=(\sigma=0)$. Note that $\l$ is generated by
global sections and the morphism given by its global sections is an
embedding on $X\setminus D$. In particular $L$ is transversal to
$\varepsilon_\cdot:(X_\cdot, D_\cdot)\to (X, D)$. Finally let
$D^L=D\resto L$, $\m_L=\m\resto L$, and $\l_L=\l\resto L$.

Now $\lambda_\cdot={\varepsilon_\cdot}\resto {L_\cdot}: (L_\cdot,
D^{L_\cdot}_\cdot)\to (L, D^L)$ is a good hyperresolution, where
$L_\cdot=X_\cdot\times_X L$. Furthermore by \quote\dbconstr\
$$
\align
 \Ox\cdot.D. &=\rpforward{\varepsilon_\cdot}\Omega^\cdot_{X_\cdot}(\log
 D_\cdot) \\
 \Ol\cdot.D^L. &=\rpforward{\lambda_\cdot}\Omega^\cdot_{L_\cdot}(\log
 D^{L_\cdot}_\cdot). \mytag\xl
\endalign
$$
\endnewnum\label\notegy


\proclaim{\num Lemma}
One has the following \dt:
$$
 \dist {\Ol p-1.D^L.}.{\Ox p.D.
\tensor^L \l_L}.{\Ol  p.D^L.\tensor \l_L}. \mytag\fsesp
$$
\endproclaim\label\fses

\proof
First assume that $X$ is smooth and $D$ is an effective normal
crossing divisor. Then one has the following commutative diagram
\cf{Esnault-Viehweg92, 2.3}:
\diagram[h=.675cm]
& & & & 0 & & 0 & & \\
& & & & \dTo & & \dTo & & \\
0 & \rTo & {\l^{-1}_L} & \rTo & \Omega^1_X\resto L & \rTo & \Omega^1_L & \rTo & 0 \\
& & \dDotsto^\alpha & & \dTo & & \dTo && \\
0 & \rTo & \Cal K & \rTo & \Omega^1_X(\log D)\resto L & \rTo_{\qquad} & \Omega^1_L(\log D^L) & \rTo & 0 \\
& & & & \dTo & & \dTo & & \\
& & & & \oplus {\ring D_j.\resto L} & \rDotsto_\beta & \oplus {\ring D_j\cap L.} && \\
& & & & \dTo & & \dTo & & \\
& & & & 0 & & 0 & & \\
\enddiagram

$L$ is transversal to $D$, so $\beta$ is an isomorphism, hence so is
$\alpha$. Taking exterior powers one obtains that for all $p$:
$$
 \ses {\Omega^{p-1}_L(\log D^L)\tensor \l_L^{-1}}.{\Omega^p_X(\log
 D)\resto {L}}.{\Omega^p_L(\log D^L)}..\mytag\smoothfses
$$

Next consider the general case. Let \letgoodhyp. 
By \quote\smoothfses\ one has the following \shes\ for all $i$:
$$
 \ses {\Omega^{p-1}_{L_i}(\log
 D^{L_i}_i)}.{\Omega^p_{X_i}(\log D_i)
 \tensor \l_{L_i}}.{\Omega^p_{L_i}(\log D^{L_i}_i)\tensor \l_{L_i}}..
$$
Since $L_i$ is the pull-back of $L$ for all $i$ these maps are
compatible with $\lambda_\cdot$, and then applying
$\rpforward{\lambda_\cdot}$ gives the required \dt.
\endproof


\section{Trace map, Gysin morphism, etc.}
The first subsection of this section is an adaptation of some parts of
\cite{Hartshorne75, II.2-3} to the logarithmic setting.  

\subsection{The trace map}
\newnum
In addition to the notation and assumptions of \quote\notegy, $X$ and
$L$ will be assumed to be smooth and $D$ an effective normal crossing
divisor throughout this subsection. Consider the following \shes,
$$
\ses \l^{-1}_L.\Omega^1_X (\log D) \resto L.\Omega^1_L (\log D^L).,
$$
and the induced natural map,
$$
 \Omega^p_L (\log D^L) \to \Omega^{p+1}_X(\log D)\resto L \otimes
 \omega_{L/X}, \mytag\hartegy
$$
where $\omega_{L/X}\simeq \l_L$ as in \cite{Hartshorne75}.

Through the rest of this section all morphisms between sheaves and
complexes are meant to be in $D(X)$ even if only sheaves are
involved. Let $\iota: L\to X$ be the embedding of $L$ into $X$. The
definition of $\iota^!$ for a finite morphism
\cite{Hartshorne66, VI 3.1, p.\ 311, p.\ 165} together with the
fundamental local isomorphism \cite{Hartshorne66, III 7.2} shows that
$$
 \iota^!\Omega_X^p (\log D)\qis \Omega^{p}_X(\log D)\resto L \otimes
 \omega_{L/X}[-1],
$$
and then the trace map for residual complexes gives:
$$ 
 \opname{Tr_\iota} : \rpforward\iota \left(\Omega^{p}_X(\log D)\resto
 L \otimes \omega_{L/X}\right)[-1] \to \Omega_X^p(\log
 D). \mytag\hartket
$$
Combining \quote\hartegy\ and \quote\hartket\ one has:
$$
\rpforward\iota \Omega^{p}_L(\log D^L)[-1]\to  \Omega^{p+1}_X(\log D). 
$$
Note, that the left hand side is supported on $L$, so the map factors
through $R\Cal H^{\vphantom{p+1}}_L\Omega^{p+2}_X(\log D)$. Also, by
the proof of \cite{Hartshorne75, II.2.2} this map is compatible with
the differential of the de Rham complex, so by taking the above for
all $p$ one has a trace map:
$$
 \tr: \rpforward\iota \Omega^{\cdot}_L(\log D^L)[-2]\to R\Cal
 H^{\phantom{\cdot}}_L\Omega^{\cdot}_X(\log D). \mytag\trace
$$
\endnewnum

\proclaim{\num Lemma \cite{Hartshorne75, II.3.1}}
$\tr$ in \quote\trace\ is a quasi-isomorphism.
\endproclaim\label\traceiso

\demo{\remnum Remark}
The proof of this lemma is taken from \cite{Hartshorne75, II.3.1} with
some small modifications and repeating it for logarithmic
differentials instead of ordinary ones. It is included for the benefit
of the reader as this constitutes an important step in the entire
proof.
\enddemo

\proof
Since $L$ is of codimension one and $\Omega^{p}_X(\log D)$ is locally
free for all $p$, there is only one non-zero local cohomology sheaf,
namely $\Cal H^1_L\Omega^{p}_X(\log D)$. Furthermore $\Cal H^1_L
\Omega^{p}_X(\log D)$ can be identified with $\Omega_X^p(\log
D)\resto{X\setminus L}/ \Omega_X^p(\log D)$, locally isomorphic to
$\Omega_X^p(\log D)[f^{-1}]/ \Omega_X^p(\log D)$ where $f$ is a local
equation of $L$ in $X$.
Therefore
$$ 
R\Cal H^{\phantom{\cdot}}_L\Omega^{\cdot}_X(\log D)\qis \Cal
H^1_L\Omega^{\cdot}_X(\log D)[-1], \mytag\trtolocal
$$
where $\Cal H^1_L\Omega^{\cdot}_X(\log D)$ denotes the complex whose
$p^{\text{th}}$ term is $\Cal H^1_L\Omega^{p}_X(\log D)$ and whose
differential is the one induced from the differential of the de Rham
complex.

Using the local isomorphism $\Cal H^1_L \Omega^{p}_X(\log D)\simeq
\Omega_X^p(\log D)[f^{-1}]/ \Omega_X^p(\log D)$ one sees easily that
$\tr$ is induced by the map
$$ 
 \Omega^{p}_L(\log D^L)[-1]\to \Cal H^1_L\Omega^{p}_X(\log D),
 \mytag\trwedge
$$
given by $\eta\mapsto\eta\wedge f^{-1}df$, where $\eta\in
\Omega^{p}_L(\log D^L)$. 

Next assume that $L=\Spec A$ and $X=\Spec B$ are affine.  Let $f\in B$
be an equation of $L$ in $X$, so $A\simeq B/(f)$. It is enough to
prove the desired quasi-isomorphism after passing to the completion
with respect to the $f$-adic topology (\cf{Hartshorne75, p.38}), so we
may assume, that $B\simeq A[[f]]$ by \cite{Hartshorne75, II.1.2}.

Based on the above discussion it will be sufficient to show that the
map,
$$
 \tau: \Omega_A^\cdot(\log D^L)[-1]\to \Omega_B^\cdot(\log
 D)[f^{-1}]/\Omega_B^\cdot(\log D),
$$
given by $\eta\mapsto\eta\wedge f^{-1}df$ for $\eta\in \Omega^{p}_L
(\log D^L)$ is a quasi-isomorphism of complexes.

Let $\gamma\in \Omega_B^{p+1}(\log D)[f^{-1}]/\Omega_B^{p+1}(\log
D)$. Then $\gamma$ can be written as $\gamma=\gamma_1+\gamma_2\wedge
df$, where $\gamma_1\in
\Omega_B^{p+1} (\log D)[f^{-1}]$ and $\gamma_2\in \Omega_B^{p}(\log
D)[f^{-1}]$, \st\ $\gamma_i=\sum_{j=0}^{k}\gamma_{ij}f^{-j}$ for
$i=1,2$ and a suitable $k\in \N$ with $\gamma_{ij}\in
\Omega_B^{p+2-i}(\log D)$. Furthermore, using the fact that $B\simeq
A[[f]]$, $\gamma_{ij}=\sum_{n=0}^{\infty}
\gamma_{ijn} f^n$ where $\gamma_{ijn}\in \Omega_A^{p+2-i}(\log
D^L)$. Hence
$$
 \gamma =\sum_{j=0}^{k} \sum_{n=0\vphantom{j}}^{\infty}
 \left(\gamma_{1jn} + \gamma_{2jn}\wedge df\right)f^{n-j}.
$$
Notice that all but a finite number of terms of this expression will
be in $\Omega_B^{p+1}(\log D)$, so one obtains that $\gamma$ can be
written uniquely in the form
$$
\gamma= \sum_{s=1}^{N}(\alpha_s + \beta_s \wedge df) f^{-s}\mytag\gam
$$
for a suitable $N\in \N$ and $\alpha_s\in \Omega_A^{p+1}(\log D^L)$
and $\beta_s\in \Omega_A^{p}(\log D^L)$.

Now $d\gamma=0$ if and only if
$$
\align
d\alpha_s &=0, \quad s=1,\dots,N \\
d\beta_1 &=0 \\
d\beta_{s+1} &=(-1)^{p+1}s\alpha_s,  \quad s=1,\dots,N.
\endalign
$$
Let $\theta=(-1)^{p+1}\sum_{s=2}^N \frac1{s-1}\beta_sf^{-s+1}\in
\Omega_A^{p}(\log D^L)$. Then $\gamma= d\theta + \beta_1\wedge
f^{-1}df$ where $d\beta_1 =0$. Hence $h^i(\tau)$ is surjective for all
$i$.

Finally if $\beta_1\wedge f^{-1}df=d\gamma'$ for some $\gamma'\in
\Omega_B^{p}(\log D)[f^{-1}]/\Omega_B^{p}(\log
D)$, then a similar expression for $\gamma'$ as the one for
$\gamma$ in \quote\gam\ shows that then there exists a $\rho\in
\Omega_A^{p-1}(\log D^L)$ \st\ $\gamma'=\rho\wedge
f^{-1}df$. Therefore $h^i(\tau)$ is also injective for all $i$.
\endproof

\subsection{Strong ampleness}
\hskip-6pt 
In this subsection the extra assumptions made in the previous
subsection are dropped, in particular $X$ is not necessarily smooth,
but \quote\notegy\ is still in effect.

\demo{\bf \num Definition}
Let $\Cal K$ be a semi-ample line bundle on $X$. Then $\Cal K$ is
called {\it strongly ample with respect to $X\setminus D$\/} if it is
ample with respect to $X\setminus D$ and \te s a proper birational
morphism, $\alpha: \tilde X\to X$, \st\ for $\tilde D=\alpha^*D$, 
$\alpha_{\tilde X\setminus\tilde D}:\tilde X\setminus\tilde D\to
X\setminus D$ is an isomorphism and there exists an effective divisor
$\tilde B$ on $\tilde X$ \st\ $\supp \tilde B= \supp \tilde D$ and
$\alpha^*\Cal K^a (\tilde B)$ is an ample line bundle for some $a>0$.

In particular if $\tau\in\coh 0.X.\Cal K^{am}.$ is a general section
for some $m\gg 0$ and $K=(\tau=0)$, then $\alpha^*K+m\tilde B$ is an
effective ample Cartier divisor supported on $\tilde D\cup \alpha^*
K$, hence $\tilde X\setminus (\tilde D\cup \alpha^* K)\simeq
X\setminus (D\cup K)$ is affine.
\enddemo\label\strongdef

\newnum Note that if $D=\emptyset$, then $\Cal K$ is strongly ample if
and only if it is ample. It is also clear that if $\Cal K$ is strongly
ample then it is also big. On the other hand, let $\pi:X\to\Bbb P^n$
be the blow up of $\Bbb P^n$ at a single point for $n\geq 2$. Let $D$
be the exceptional divisor of $\pi$. Then $\pi^*\ring X.(1)$ is
semi-ample and big, but not strongly ample with respect to $X\setminus
D$.
\endnewnum\label\nonstrong

It will be very important in \S 4 that this property is inherited by
restrictions to $L$:

\proclaim{\num Lemma}
If $\Cal K$ is strongly ample with respect to $X\setminus D$, then
$\Cal K_L=\Cal K\resto L$ is strongly ample with respect to
$L\setminus D^L$.
\endproclaim\label\strongres

\proof
$\Cal K$ is ample with respect to $X\setminus D$, so $\Cal K_L$ is
ample with respect to $(X\setminus D)\cap L=L\setminus D^L$. Let
$\alpha: \tilde X\to X$ be a proper birational morphism, $\tilde
D=\alpha^* D$, and $\tilde B$ an effective divisor on $\tilde X$ \st\
$\supp \tilde B=\supp \tilde D$ and $\Cal\alpha^*\Cal K^a(\tilde B)$
is ample for some $a>0$.  

Let $\tilde L$ be the proper transform of $L$ on $\tilde X$,
$\alpha_L=\alpha\resto{\tilde L}$ and $\tilde B^L=\tilde B\resto
L$. It is easy to see that $\Cal K_L$ and $\alpha_L:\tilde L\to L$
satisfies the requirements of the definition
\quote\strongdef.
\endproof

The following lemma gives important examples for strongly ample line
bundles.

\proclaim{\num Lemma}
Assume that \te s an effective $\Q$-Cartier divisor $B$, \st\ $\supp
B=\supp D$ and, in addition to \quote\notegy, one of the following
holds:
\myroster
\myitem
$\l$ is ample, or
\rostlabel\affineegy
\myitem
$B$ is nef.
\rostlabel\affineket
\endmyroster
Then $\m$ is strongly ample with respect to $X\setminus D$.
\endproclaim\label\strongex

\proof
It is enough to prove that $\l=\m^N$ is strongly ample with respect to
$X\setminus D$.

If $\l$ is ample and $B$ is $\Q$-Cartier, then $\l^a(bB)$ is ample for
some $a, b>0$.

In case \quote\affineket, let $\phi:X\to Z$ be the morphism given by
the global sections of $\l$ and $\Cal A$ an ample line bundle on $Z$,
\st\ $\l=\phi^*\Cal A$. Further let $\alpha:\tilde X\to X$ be the
blowing up of $X$ along the exceptional set of $\phi$ and $\tilde
D=\alpha^*D$. Note that the exceptional set of
$\tilde\phi=\phi\circ\alpha$ is a Cartier divisor with support
contained in $\supp \tilde D=\supp \alpha^* B$.

Now there exists a $\tilde \phi$-exceptional divisor $E$ on $\tilde
X$, \st\ $\tilde \phi^*\Cal A^a(-E)=\alpha^*\l^a(-E)$ is ample for
some $a>0$.  Since $E$ is $\tilde\phi$-exceptional, $\supp
E\subseteq\supp \tilde D=\supp \alpha^* B$, so for some $b>0$, $\tilde
B=b\alpha^* B-E$ is effective and $\supp \tilde B=\supp\tilde D$.
Finally $\alpha^*\l^a(\tilde B)$ is ample, since $B$ is nef.
\endproof

\subsection{Gysin morphism}
\newnum
Using the notation and assumptions of \quote\notegy\ further assume
that $\m$ is strongly ample with respect to $X\setminus D$.
\endnewnum\label\noteket

\newnum
Again, $\varepsilon_\cdot :(X_\cdot, D_\cdot)\to (X,D)$ denotes a good
hyperresolution. Applying
\quote\traceiso\ for $\iota^i: L_i\hookrightarrow X_i$
one obtains the following natural quasi-isomorphism
$$
 \rpforward\iota^i \Omega^{\cdot}_{L_i} (\log
 D_i^{L_i})[-2] @>{\phantom{_{q}}\qis}>> R\Cal
 H^{\phantom{\cdot}}_{L_i}\Omega^{\cdot}_{X_i}(\log
 D_i).
$$
By \quote\xl\ this implies that there exists a quasi-isomorphism
$$
 \rpforward\iota \underline{\Omega}^{\cdot}_{L} (\log
 D^{L})[-2] @>{\phantom{_{q}}\qis}>> R\Cal
 H^{\phantom{\cdot}}_{L}\underline{\Omega}^{\cdot}_{X}(\log
 D). \mytag\Rgysin
$$

Let $j:X\setminus D\to X$ and $j^L: L\setminus D^L\to L\to X$ be the
inclusion maps. Then by \quote\cc\ \quote\Rgysin\ gives a
quasi-isomorphism
$$
 R{j^L_*} \C_{L\setminus D^L}[-2] @>{\phantom{_{q}}\qis}>> 
R\Cal H^{\phantom{\cdot}}_{L}\rpforward{j}\C_{X\setminus D}.
\qis
R\Cal H^{\phantom{\cdot}}_{L\setminus D^L}\C_{X\setminus D}.
$$
Applying $R\Gamma$ to both sides one obtains a quasi-isomorphism
$$
 R\Gamma \C_{L\setminus D^L} [-2]@>{\phantom{_{q}}\qis}>>
 R\Gamma_{L\setminus D^L}\C_{X\setminus D}.
$$
In particular
$$
\coh i-2.L\setminus D^L. \C.@>{ \simeq }>>
\loccoh i.{L\setminus D^L}.X\setminus D. \C.\mytag\gysin
$$ is an isomorphism for all $i$.

On the other hand, 
$$ 
 \dist {R\Gamma_{L\setminus D^L} \C_{X\setminus
D}}.R\Gamma\C_{X\setminus D}.{R\Gamma\C_{X\setminus (D\cup L)}}
.
$$
forms a \dt.  By \quote\strongdef\ $X\setminus (D\cup L)$ is
affine, so $\coh j. X\setminus (D\cup L).\C.=0$ for $j>\dim X$
\cite{Hartshorne75, II.4.6}, \cite{Goresky-MacPherson83},
\cite{GNPP88, III.3.1(i)}. Hence
$$
\loccoh i. {L\setminus D^L}. X\setminus D.\C.
\to \coh i. X\setminus D.\C.
$$
is an isomorphism for $i>\dim X+1$ and surjective for $i=\dim
X+1$. Combinig this with \quote\gysin\ one obtains that
$$
\coh i-2.L\setminus D^L. \C. \to
\coh i. X\setminus D.\C.
$$
is an isomorphism for $i> \dim X+1$ and surjective for $i=\dim
X+1$. Furthermore by the construction of these maps it is clear that
they respect the Hodge decomposition \quote\hodge. Therefore
$$
G: \hypcoh q-1.L.{\Ol p-1. D^L.}. \to \hypcoh q.X.{\Ox
 p. D.}. \mytag\indstep
$$
is an isomorphism for $p+q> \dim X+1$ and surjective for $p+q=\dim
X+1$.
\endnewnum
\endsubsection

\section{The logarithmic Kodaira-Akizuki-Nakano Vanishing Theorem}

\proclaim{\num Theorem}
With the notation and assumptions of \quote\notegy\ and
\quote\noteket,
$$
\hypcoh q.X.{\Ox p.D.\otimes \l}.=0\quad \text{ for } p+q>\dim X.
$$
\endproclaim\label\kanvery

\proof
Tensoring the \shes,
$$
\ses {\ring X.}.\l.\l_L.,
$$
by $\Ox p.D.$ leads to the \dt,
$$
\dist {\Ox p.D.}.{\Ox p.D.}\otimes\l.{\Ox p.D.}\otimes^L\l_L.
$$
and the corresponding long exact hypercohomology sequence:
$$ 
 \hypcoh q-1.L.{\Ox p.D.\otimes^L\l_L}. @>\partial>> \hypcoh q.X.{\Ox
 p.D.}. @>>> \hypcoh q.X.{\Ox p.D.\otimes\l}.. \mytag\psidef
$$

On the other hand,
\quote\fses\ gives the \dt,
$$
\dist {\Ol p-1.D^L.}.{\Ox p.D. 
\tensor^L \l_L}.{\Ol p.D^L.\tensor \l_L}.,
$$
and in turn the long exact hypercohomology sequence:
$$ 
\hypcoh q-1.L.{\Ol p-1.D^L.}. \to\hypcoh q-1.L.{\Ox
p.D.\otimes^L\l_L}. \to\hypcoh q-1.L.{\Ol p.D^L.\otimes\l_L}..
$$ 

Now by induction and \quote\strongres\ we may assume that $\hypcoh
q-1.L.{\Ol p.D^L.\otimes\l_L}.=0$. Hence
$$
 \phi: \hypcoh q-1.L.{\Ol p-1.D^L.}. \to\hypcoh q-1.L.{\Ox
 p.D.\otimes^L\l_L}.
$$
is an isomorphism for $p+q>\dim X +1$ and surjective for $p+q=\dim
X+1$. Furthermore $\phi$ is induced by the map ${\Ol p-1.D^L.}\to {\Ox
p.D.\otimes\l_L}$ which locally is given by $\eta\mapsto \eta\wedge
df\tensor f^{-1}$ where $f$ is a local equation of $L$ in $X$. So
$\phi$ is defined the same way as the Gysin map was, hence the
following diagram is commutative, where $G$ is from
\quote\indstep\ and $\partial$ is from
\quote\psidef.

$$
\diagram[h=.75cm]
\hypcoh q-1.L.{\Ol p-1.D^L.}. && \\
\dTo_{\phi}  & \rdTo^{G} & \\
{\hypcoh q-1.L.{\Ox p.D.\otimes^L\l_L}.} 
& \rTo_{\quad\partial\quad} &
{\hypcoh q.X.{\Ox p.D.}.}  \\
\enddiagram
$$

Now $G$ and $\phi$ are isomorphisms for $p+q>\dim X +1$ and surjective
for $p+q=\dim X+1$, so the same holds for $\partial$. However, then
\quote\psidef\ implies that
\medskip
$
\hfill {\hypcoh q.X.{\Ox p. D.}\tensor\l.=0\quad \text{ for } p+q>\dim
 X}. \hfill\square
$
\enddemo

To obtain the statement in the general case one uses the usual
covering trick:

\proclaim{\num  Corollary}
$
\hypcoh q.X.{\Ox p. D.}\tensor\m.=0\quad \text{ for } p+q>\dim
X.$
\endproclaim\label\cover

\proof
Let $$
\pi: \tilde X=\spec X.\bigoplus_{i=0}^{N-1}\m^{-i}\to X.
$$ 
be the cover obtained by taking the $N^{\text{th}}$-root of $L$.  Now
the trace map of $\pi$ provides a left inverse to the natural map
\cf{GNPP88, p.151}, \cite{Esnault-Viehweg92, 3.22}:
$$
\Ox p.D.\to R\pi_*\underline{\Omega}_{\tilde X}^p(\log \pi^*D).
$$
Applying $\Bbb H^q$ and using \quote\kanvery\ on $\tilde X$ proves the
statement.
\endproof

Therefore we have (no longer using \quote\notegy):

\proclaim{\num Theorem}
Let $X$ be a projective variety and $D$ an effective 
divisor on $X$.
Let $\l$ be a semi-ample line bundle on $X$ that is strongly ample
with respect to $X\setminus D$. Then for $p+q>n$,
$$
 \hypcoh q.X.{\Ox p. D.}\tensor\l. =0.
$$
\endproclaim\label\logkan

\proclaim{\num Corollary=\gr\ Theorem}
Let $Y$ be a smooth complex projective variety and $\m$ a
semi-ample and big line bundle on $Y$. Then $$\coh
i.Y.\omega_Y\tensor\m.=0\qquad\text{for }i>0.$$
\endproclaim

\proof
First assume that $\m$ is generated by global sections. Let $\phi:Y\to
X$ be the morphism given by the global sections of $\m$. Then \te s an
ample line bundle, $\l$, on $X$ \st\ $\m=\phi^*\l$, so
by \quote\topomega, \quote\strongex, and
 \quote\logkan
$$
\coh i.Y.\omega_Y\tensor\m.\simeq 
\hypcoh i.X.R\phi_*\omega_Y\tensor\l.\simeq
\hypcoh i.X.{\Ox n. \emptyset.}\tensor\l. =0.
$$
The generic case is now proved by the usual covering trick cf.\
\quote\cover.
\endproof

\demo{\remnum Remark}
This shows that \quote\logkan\ is indeed a common generalization of
\quote\kodaira, \quote\kanlabel, \quote\gr, \quote\navarro\ and
\quote\ev.  
\enddemo

We also have a local version:

\proclaim{\num Theorem}
Let $\psi:X\to Z$ be a projective morphism and $D$ an effective $\Bbb
Q$-Cartier divisor on $X$.  Let $\l$ be a $\psi$-ample line bundle on
$X$. Then for $p+q>n$,
$$
R^q\psi_*(\Ox p.D.\tensor\l)=0.
$$
\endproclaim

\proof
The statement is local, so we may assume that $Z$ is projective. Let
$\m$ be an ample line bundle on $Z$, \st\ for all $p,q$,
$R^q_{\vphantom{X}}\psi_*(\Ox p.D.\tensor\l)\tensor\m$ is
generated by global sections and have no higher cohomology. This can
be done because $\Ox p.D.$ is bounded and has coherent cohomology
sheaves.  Furthermore choose $\m$ in such a way, that
$\l\tensor\psi^*\m$ be ample on $X$. Then 
by the Leray spectral sequence,
\quote\strongex\ and \quote\logkan,
$$
\coh 0.Z.{R^q\psi_*(\Ox p.D.\tensor\l)\tensor\m}. =
 \hypcoh q.X.{\Ox p. D.}\tensor\l\tensor\psi^*\m. =0.
$$
Since $R^q\psi_*(\Ox p.D.\tensor\l)\tensor\m$ is generated by global
sections, this proves the statement.
\endproof

Finally, this gives a bound on the range of degrees where $\Ox p.D.$
can have non-zero cohomology sheaves. 

\proclaim{\num Corollary}
Let $X$ be a projective variety and $D$ an effective $\Bbb Q$-Cartier
 divisor on $X$.

Then $h^q_{\vphantom{X}}(\Ox p.D.)=0$ for $q>n-p$ or $0>q$.
\endproclaim

\proof
Let $\psi=\id_X:X\to X$ and $\m=\ring X.$. The second inequality
is simply \quote\dbconstr.
\endproof

Regarding the case of $\coh q.Y.\Omega^p_Y\tensor\m.$ for $p<n$,
Ramanujam has already noticed that if $\m$ is only semi-ample (or
even generated by global sections) and big, than vanishing does not
necessarily hold \cite{Ramanujam72}.
However, since globally generated and big line bundles are pull-backs
of (very) ample ones, \quote\logkan\ can be considered as a
substitute. Later applications will show that it can actually be used
for this purpose.

\section{Relative complexes}
Let $f:X\rightarrow C$ be a morphism such that $C$ is a smooth complex
curve. Let $\Delta\subseteq C$ be a finite set and $D=f^*\Delta$. Let
\letgoodhyp, and consider the map $f_i=f\circ\varepsilon_i:
X_i\to C$.

The goal is to construct a complex whose cohomological properties
resemble those of $\Omega_{X/C}^{p}$ in the smooth case.

Taking the wedge product induces a map,
$$
\Omega^p_{X_i}(\log D_i)\tensor f_i^*\Oc\Delta. 
\to \Omega^{p+1}_{X_i}(\log D_i).
$$
This is obviously compatible with $\varepsilon_\cdot$, so it gives a
\moc:
$$ 
\wedp : \Ox p.D.\otimes f^*\Oc\Delta.\to \Ox p+1.D..
$$
It is also easy to see that this is independent of the actual
hyperresolution used \cf{Kov\'acs96, p.375}.  Hence $\wedp$ is a
well-defined natural map in $D(X)$.

Choose a representative, $K_{p}\in\Cal Obj(C(X))$, of $\Ox p.D.$ for
all $p$ \st\ $\wedp$ is represented by morphisms $K_{p}\to
K_{p+1}\text{ in }\Cal Mor(C(X))$.  By abuse of notation this will
also be denoted by $\wedp$.  Let $\wedp'=\wedp\otimes
id_{f^*\Oc\Delta.}\in\hom_{C(X)}(\Oxx p.D.\otimes
f^*\Oc\Delta., \Oxx p+1.D.\otimes f^*\Oc\Delta.)$. Since $\Oc\Delta.$
is a line bundle, $\wedp\circ\wedge'_{p-1}=0$.  Let
$M_{r}=0\in\Cal Obj(C(X))$, $w''_{r}=0\in
{\operatorname{Hom}}_{C(X)} (\Oxx r.D.\otimes f^*\Oc\Delta.,
M_{r}\otimes f^*\Oc\Delta.)$ and $w'_{r}=0\in
{\operatorname {Hom}}_{C(X)} (M_{r}\otimes
f^*\Oc\Delta., \Oxx r+1.D.)$ for $r\geq n$.  Assume that $p<n$ and for
every $q>p$, $M_{q}\in\Cal Obj(C(X))$ is
defined. Assume further that there are morphisms of complexes,
$$
w''_{q}: \Oxx q.D.\otimes f^*\Oc\Delta. \rightarrow
M_{q}\otimes f^*\Oc\Delta.
\qquad\text{and}\qquad 
w'_{q}: M_{q}\otimes f^*\Oc\Delta. \rightarrow \Oxx q+1.D., 
$$ 
such that 
$$
\wedge_{q}=w'_{q}\circ w''_{q}\qquad\text{and}\qquad
w''_{q}\circ\wedge'_{q-1}=0.  
$$ 
Let 
$$
w_{q}=w''_{q}\otimes id_{f^*\Oc\Delta.^{-1}} : \Oxx r.D.\to
M_{r}
$$ 
and
$$
M_{p}=M(w_{p+1})[-1]\otimes f^*\Oc\Delta.^{-1}\in\Cal Obj(C(X)),
$$
i.e.,
$$ 
M_{p}^{m}\otimes f^*\Oc\Delta. = {\Oxx p+1.D.^m}\oplus
M_{p+1}^{m-1}
$$ 
and 
$$
d^{m}_{M_p\otimes f^*\Oc\Delta.}=
\left( \matrix
\ d_{\Oxx p+1.D.}^{m} & 0 \\
-w_{p+1}^{m} & -d_{M_{p+1}}^{ m-1}  
\endmatrix \right). 
$$ 
Also let 
$$ 
w''_{p} = \left( \matrix \wedp \\ 0 \endmatrix \right):
\Oxx p.D.\otimes f^*\Oc\Delta.  \rightarrow M_{p}\otimes f^*\Oc\Delta. 
$$ 
and 
$$
w'_{p} = (id_{\Oxx p+1.D.},0): M_{p}\otimes
f^*\Oc\Delta. \rightarrow \Oxx p+1.D..
$$ 

\noindent $w'_{p}$ is a \moc\ by the definition of the mapping cone
and $w''_{p}$ is a \moc\ because $w_{p+1}\circ\wedge_{p}=0$. It is
also obvious that $\wedp = w'_{p}\circ w''_{p}$ and
$w''_{p}\circ\wedge'_{p-1}=0$.

\noindent Also, by their definition, the equivalence classes of $w_p$,
$w'_p$ and $w''_p$ in $D(X)$ are independent of the hyperresolution
chosen.  From now on these symbols will denote their equivalence
classes in $D(X)$. A map will mean an element of $\Cal Mor(D(X))$, so
it is possibly not represented by an actual morphism of complexes
between two arbitrary representatives of the respective objects.

\proclaim{\num Theorem-Definition}
Let $f:X\rightarrow C$ be a morphism between complex varieties such
that ${\operatorname {dim}}X=n$ and $C$ is a smooth curve. Let
$\Delta\subseteq C$ be a finite set and $D=f^*\Delta$. For every
nonnegative integer p there exists a complex $\,
\Oxcp\in\Cal Obj(D(X))$ with the following properties.

\myroster
\myitem
The natural map $\wedge_{p}$ factors through $\, \Oxcp\otimes
f^*\Oc\Delta.$, i.e.,\ there exist maps:
$$
\align
w''_{p} &:\Oxp\otimes f^*\Oc\Delta.\rightarrow\Oxcp\otimes
f^*\Oc\Delta. \qquad\text{and}\\ 
w'_{p} &:\Oxcp\otimes f^*\Oc\Delta. \rightarrow \Ox p+1.D. 
\endalign
$$ 
such that
$\wedge_{p}=w'_{p}\circ w''_{p}$.  

\myitem
If $w_{p}=w''_{p}\otimes
id_{f^*\Oc\Delta.^{-1}} : \Oxp\rightarrow\Oxcp$, then 
$$ 
 \Oxcp\otimes f^*\Oc\Delta. @>~w'_{p~}>> \Ox {p+1}.D. @>w_{p+1}>>
 \Oxc{p+1}. @>{~+1~}>>
$$ 
is a distinguished triangle in
$D(X)$.  
\rostlabel\Odist

\myitem
$w_p$ is functorial, i.e.,\ if $\phi :Y\rightarrow X$ is a C-morphism,
then there are natural maps in $D(X)$ forming a commutative diagram:

$$
\diagram[h=.675cm]
\ \ \Oxp & \rTo & \ \ \ \Oxcp \\ 
\dTo && \dTo \\ 
R\phi_{*}\Oy {p}.\phi^*D. & \rTo_{\qquad} &  R\phi_{*}\Oyc
{p}.\phi^*D. 
\enddiagram
$$ \hfill
\rostlabel\pushsquare

\myitem
If $f$ is smooth over $C\setminus\Delta$, then $\Oxcp\qis
\Omega_{X/C}^{p}(\log D)$
\rostlabel\Osmooth

\myitem
$\Oxc r.=0$ for $r\geq n$ and if f is proper, then $\, \Oxcp\in\Cal
Obj(D^{b}_{coh}(X))$ for every p.
\rostlabel\largexc


\endmyroster
\endproclaim 

\proof Let $\Oxcp\qis [M_{p}] \in \Cal Obj(D(X))$.
Then (\the\sect.\the\thm.1), (\the\sect.\the\thm.2) and the first part
of (\the\sect.\the\thm.5) follows from the discussion above.  Using
(\the\sect.\the\thm.2), the first part of (\the\sect.\the\thm.5) and
descending induction on $p$, (\the\sect.\the\thm.3),
(\the\sect.\the\thm.4) and the rest of (\the\sect.\the\thm.5) follows
from \quote\pushmap\ and \quote\dbsmooth.  
\endproof

Note that the combination of \quote\pushmap, \quote\zeromap, and
\quote\pushsquare\ implies that 
if $\phi :Y\rightarrow X$ is a C-morphism, then there are natural maps
in $D(X)$ forming a commutative diagram:
$$
\diagram[h=.675cm]
\ring X. & \rTo & \ \ \Ox 0.D. & \rTo & \ \ \ \Oxc 0. \\ 
\dTo & & \dTo & & \dTo \\ 
R\phi_{*}\ring Y. & \rTo_{\qquad\qquad} &
R\phi_{*}\Oy {0}.\phi^*D.  & \rTo_{\qquad\qquad} & R\phi_{*}\Oyc
{0}.\phi^*D. 
\enddiagram \mytag\cdforvan
$$

\section{More vanishing theorems}
\proclaim{\num Theorem}
Let $X$ be a projective variety of dimension $n$ and $f:X\rightarrow
C$ a morphism to a smooth complex curve. Let $\Delta\subseteq C$ be a
finite set and $D=f^*\Delta$.  Let $\l$ be a line bundle on $X$ \st\
$\l$ and $\l\tensor f^*\Oc\Delta.^{-(n-1)}$ are semi-ample and
ample with respect to $X\setminus D$.  Then
$$
\hypcoh n.X.{\Oxc 0.\tensor\l\tensor f^*\Oc\Delta.}.=0.
$$
\endproclaim\label\lastquote

\proof
Let $\l_p=\l\tensor f^*\Oc\Delta.^{-(p-1)}$ for $p=0,\dots,n$. By
assumption, $\l_p$ is semi-ample and ample with respect to $X\setminus
D$ for $1\leq p \leq n$ since either $f^*\Oc\Delta.$ or
$f^*\Oc\Delta.^{-1}$ is semi-positive. In fact $\l_p$ is strongly
ample with respect to $X\setminus D$ for $1\leq p \leq n$ by
\quote\strongex\ since $D$ is an effective, nef $\Bbb Q$-Cartier
divisor.
Twisting  \quote\Odist\ by $\l_p$ yields the following \dt:
$$
\Oxc p-1.\tensor\l_{p-1} @>>> \Ox p.D.\tensor\l_{p} @>>>
\Oxc{p}.\tensor\l_{p} @>{~+1~}>>. 
$$ 

By \quote\logkan\ $\hypcoh n-(p-1).X.{\Ox p. D.}\tensor\l_p. =0$, so
the map 
$$ 
\hypcoh n-p.X.{\Oxc p.}\tensor\l_p. \to 
\hypcoh n-(p-1).X.{\Oxc p-1.}\tensor\l_{p-1}. 
$$ 
is surjective for all $1\leq p \leq n$. Observe that these maps form a
chain as $p$ runs through $p=n, n-1, \dots,1$. So the composite map
$$ 
\hypcoh 0.X.{\Oxc n.}\tensor\l_n. \to 
\hypcoh n.X.{\Oxc 0.}\tensor\l_{0}. 
$$ 
is also surjective. However $\Oxc n.=0$ by construction (cf.\
\quote\largexc). Therefore the statement follows.
\endproof

\proclaim{\num Lemma}
Let $\phi:Y\to X$ be a proper generically finite map of varieties of
dimension $n$. Let $\Cal F$ be a coherent sheaf on $Y$. Then the
natural map
$
\coh n.X.\phi_*\Cal F.\to \coh n.Y.\Cal F.
$
is surjective.
\endproclaim\label\topcohequal

\proof
Let $x\in X$ and let $d(x)=\dim Y_x$, the dimension of the fiber of
$\phi$ over $x$. Now $(R^j\phi_*\Cal F)_x=0$ for $j>d(x)$, so $\supp
R^j\phi_*\Cal F\subseteq X_j=\{ x\in X \mid d(x)\leq j\}$.  Clearly,
$X=X_0\cup X_1\cup\dots\cup X_{n-1}$ and for all $j>0$, the dimension
of $\phi^{-1}(X_j)$ is at most $n-1$, so $\dim X_j +j\leq n-1$. Hence
$\dim\supp R^j\phi_*\Cal F\leq n-j-1$ for $j>0$. 
Therefore $\coh i.X.R^j\phi_*\Cal F.=0$ for $j>0, i+j\geq n$. Finally
this implies that in the Leray spectral sequence $\coh
i.X.R^j\phi_*\Cal F.\Rightarrow \coh i+j.Y.\Cal F.$ the only non-zero
term for $i+j\geq n$ is $\coh n.X.\phi_*\Cal F.$.
\endproof

\proclaim{\num Lemma}
Let $\phi:Y\to X$ be a proper generically finite map of normal
varieties of dimension $n$. Let $\l$ be a line bundle on $X$.
\myroster
\myitem
If $\phi$ is birational, then the natural map
$
\coh n.X.\l.\to \coh n.Y.\phi^*\l.
$
is surjective.
\rostlabel\hypstepegy

\myitem
If $X$ is projective and has rational singularities, then the natural
map $
\coh n.X.\l.\to \coh n.Y.\phi^*\l.
$ is injective.
\rostlabel\hypstepket

\endmyroster
\endproclaim\label\hypstep

\proof
If $\phi$ is birational, then $\ring X.\simeq \phi_*\ring Y.$, so
$\l\simeq\phi_*\phi^*\l$, hence \quote\topcohequal\ implies
\quote\hypstepegy.
If $X$ is projective and has rational singularities, then the natural
map $\ring X.\to R\phi_*\ring Y.$ has a left inverse by
\cite{Kov\'acs00b, Theorem 2}. Hence
$
 \coh n.X.\l.\to \hypcoh n.X. {R\phi_*\ring Y.\tensor\l}.\simeq
 \coh n.Y.\phi^*\l.
$ is injective.
\endproof

\def\myhskip{\empty
}
\def\mysub{\rInto}
\def\mysup{\lInto}

\proclaim{\num Theorem}\label\momore Let $X$ be a projective variety
of dimension $n$ and $f:X\rightarrow C$ a morphism to a smooth proper
curve. Let $\Delta\subseteq C$ be a finite set and
$D=f^*\Delta$. Assume that \te s a smooth projective variety, $Y$, and
a proper generically finite map, $\phi: Y\to X$, \st\
$h\resto{Y\setminus B}:{Y\setminus B}\to C\setminus \Delta$ is smooth,
where $h=f\circ \phi$ and $B=h^*\Delta$.  Let $\tilde C$ be smooth
proper curve, $\sigma:\tilde C\to C$ a finite cover, unramified over
$C\setminus\Delta$. Assume that for $\tilde\Delta=
(\sigma^*\Delta)_{\text{red}}$, $\omega_{\tilde C}(\tilde \Delta)\subset
\sigma^*\omega_C(\Delta)$.  Let $\tilde X$ be the normalization of
$\tilde C\times_C X$ and $\tilde Y\to \tilde C\times_C Y$ a resolution
of singularities \st\ it is an isomorphism over $\tilde C\setminus
\sigma^*\Delta${\rm :} 
\diagram [size=.75cm] 
\tilde\pi^*B=\tilde\phi^*\pi^*D
\myhskip &
\myhskip \mysub\myhskip & \myhskip \tilde Y \myhskip &&&&& \\
\dTo &&\dTo^{\tilde\phi} &
\rdTo[leftshortfall=.425cm][rightshortfall=.25cm](5.225,1.9)^{\tilde\pi}
&&&& \\ 
\pi^*D=\tilde f^*\sigma^*\Delta 
\myhskip&\myhskip \mysub\myhskip &\myhskip {\tilde X}\myhskip
&&&& & Y &\myhskip \mysup\myhskip & B=\phi^*D \\ 
&\hskip-.1cm\rdTo[rightshortfall=.35cm](3.22875,4.1)
&&\rdTo[leftshortfall=.425cm][rightshortfall=.25cm](5.225,1.9)^\pi
\hskip-.1cm\rdTo[rightshortfall=.35cm](3.22875,4.1)_{\tilde f\ }
\rdTo[rightshortfall=.5cm](3,2) &&&& \dTo_\phi & & \dTo \\ 
&&&& & \tilde
C\times_C X & \rTo & X  &\myhskip \mysup\myhskip & D=f^*\Delta
\\ &&&& & \dTo & & \dTo_f& & \dTo \\ 
&& &
\sigma^*\Delta \myhskip &\myhskip \mysub\myhskip & \myhskip \tilde C &
\rTo^\sigma & C &\myhskip \mysup\myhskip & \Delta 
\enddiagram
Assume that \te s a line bundle $\l$ on $X$ \st\ $\pi^*\l$ contains a
line bundle $\tilde \l$ \st\ $\tilde\l$ and $\tilde\l\tensor\tilde
f^*\Oc\tilde\Delta.^{-(n-1)}$ are semi-ample and ample with respect to
$\tilde X\setminus \pi^*D$.

\myroster
\myitem
If $\phi$ is birational, then
$$
\coh n.Y.{\phi^*\l\tensor h^*\Oc\Delta.}.=0.
$$
\rostlabel\mainegy

\myitem
If $X$ has rational singularities, then
$$
\coh n.X.{\l\tensor f^*\Oc\Delta.}.=0.
$$
\rostlabel\mainket

\endmyroster
\endproclaim

\proof 
Let $\Cal M=\l\tensor f^*\Oc\Delta.$ and $\tilde\Cal M=\tilde\l\tensor
\tilde f^*\omega_{\tilde C}(\tilde \Delta)\subseteq \pi^*\Cal M$. By
\quote\cdforvan\ one has the following commutative diagram:
\diagram 
\coh n.\tilde X.{\tilde\Cal M}. & & \rTo & & \hypcoh n.\tilde
X.{\underline\Omega_{\tilde X/\tilde C}^0(\log \pi^* D)\tensor\tilde\Cal
M}. \\ \dTo^{\tilde\alpha} & \hskip-1cm \rdTo(7,2)^{\tilde\beta} & & &
\dTo \\ \coh n.\tilde Y.{\tilde\phi^*\tilde\Cal
M}. &\rTo^{\quad\simeq\quad}& \hypcoh n.\tilde X.{R\phi_*\ring \tilde
Y.\tensor\tilde \Cal M}. & \rTo_{\quad\eta\quad} & \hypcoh n.\tilde
X.{R\phi_*\underline\Omega_{\tilde Y/\tilde C}^0(\log
\tilde\pi^*B)\tensor\Cal M}.
\enddiagram 

 $\hypcoh n.\tilde X.{\underline\Omega_{\tilde X/\tilde C}^0(\log
\pi^* D)\tensor\tilde\Cal M}.=0$ by \quote\lastquote, so $\tilde\beta$
is the zero map. Furthermore $\eta$ is an isomorphism by
\quote\Osmooth, hence $\tilde\alpha$ is the zero map as well.

Now consider the following commutative diagram:
\diagram
 \coh n.\tilde X.{\tilde\Cal M}. & \rTo^{\gamma} & 
\coh n.\tilde X.{\pi^*\Cal M}. & \lTo^\delta & \coh n. X.{\Cal M}.  \\
\dTo^{\tilde\alpha} & 
& \dTo^{\beta} & & \dTo^\alpha \\ 
 \coh n.\tilde
 Y.{\tilde\phi^*\tilde\Cal M}. &\rTo^{\quad \tilde\gamma\quad }& \coh
 n.\tilde Y.{\tilde\phi^*\pi^*\Cal M}. & \lTo^{\quad \tilde\delta\quad } &
 \coh n. Y.{\phi^*\Cal M}.  \\
\enddiagram

The cokernel of the inclusion $\tilde\Cal M\subseteq \pi^*\Cal M$ is
supported in codimension 1, so $\gamma$ and $\tilde\gamma$ are
surjective. Therefore $\beta$ is the zero map, since so is $\tilde
\alpha$.

Now $\delta$ is injective, because $\pi$ is finite and $X$ is normal
and $\tilde \delta$ is injective by \quote\hypstepket. Hence $\alpha$
is the zero map.

If $\phi$ is birational, then $\alpha$ is also surjective by
\quote\hypstepegy, so \quote\mainegy\ follows.

If $X$ has rational singularities, then $\alpha$ is also injective by
\quote\hypstepket, so \quote\mainket\ follows.
\endproof

\section{Arakelov-Parshin boundedness}

\demo{\bf \num Definition}
A morphism, $h:Y\to C$ is called {\it isotrivial\/} if all but
finitely many fibers of $h$ are isomorphic to a fixed variety.
Similarly $h$ is called {\it birationally isotrivial\/} if all but
finitely many fibers of $h$ are birational to a fixed variety.
\enddemo

\demo{\bf\num Definition}\cite{Esnault-Viehweg90}
Let $F$ be a normal Gorenstein variety with rational singularities,
$\l$ a line bundle on $F$ and $\Gamma$ an effective divisor \st\
$\l=\ring F.(\Gamma)$. Let
\roster
\item
$\Cal C(\Gamma, N)=\opname{coker}\left(\tau_*\omega_{\tilde
F}\left(-\left[\frac{\tilde \Gamma}N\right]\right)\to \omega_F\right)$
where $\tau:\tilde F\to F$ is a resolution of singularities \st\
$\tilde \Gamma=\tau^*\Gamma$ is a normal crossing divisor.
\item
$e(\Gamma)=\min\{N\in \N_+ | \Cal C(\Gamma, N)=0\}$
\item
$e(\l)=\sup\{e(\Gamma) | \exists\lambda\in\coh 0.F.\l.\text{ \st\ }
\Gamma=(\lambda=0)\}$
\endroster
$e(\l)$ will be called the {\it Esnault-Viehweg threshold\/} of
$\l$. For properties of $e(\l)$ the reader should consult
\cite{Viehweg95, \S 5.3-4}
\enddemo\label\defofe

\demo{\bf\num Notation and Assumptions}
Let $C$ be a smooth projective curve of genus $g$, $\Delta\subseteq C$
a finite set of points, regarded as a (reduced) divisor. Let
$\delta=\#\Delta$, the number of points in $\Delta$.  Let $X$ be an
irreducible projective variety of dimension $n$ with rational
Gorenstein singularities and $f:X\to C$ a morphism. Let $D=f^*\Delta$.

Assume that $X_t$ has rational Gorenstein singularities for $t\in
\C\setminus \Delta$ and that \te s a simultaneous resolution of
$X\setminus D \to C\setminus \Delta$, i.e., \te s a smooth projective
variety $Y$ and a birational 
morphism $\phi:Y\to X$, \st\ $Y\setminus \phi^*D\to C\setminus\Delta$
is smooth. Let $h=f\circ \phi$, $B=\phi^*D$ and let $Y_{\text{gen}}$
denote the general fiber of $h$. By blowing up $Y$ along a subvariety
of $B$ one may assume that $B$ is a (not necessarily reduced) normal
crossing divisor.

$r(m)$ will denote the rank of $f_*\omega_{X/C}^m$. This is equal to
the $m$-th plurigenus of the general fiber of $f$,
$P_m(X_{\text{gen}})$.  $e(m)$ will denote the Esnault-Viehweg
threshold of $\omega_{X_{\text{gen}}}^m$.  If
$\omega_{X_{\text{gen}}}$ is ample, then $e(m)\leq m^n
K_{X_{\text{gen}}}^n +1$ for $m\gg 0$ by \cite{Viehweg95, 5.12}.
\enddemo\label\assumptions

\demo{\remnum Remark} 
 $X_t$ has only rational singularities for $t\in C\setminus \Delta$,
 so the same holds for $X\setminus D$.  It is conjectured that a
 variety with only rational singularities admits a compactification
 with only rational singularities.  Furthermore, if that conjecture
 holds then the Gorenstein assumption could be avoided as well with a
 little care. Hence the assumption on the singularities of $X$ is
 conjecturally superfluous.
\enddemo


The following lemma gives an effective measure of the positivity of
$f_*\omega_{X/C}^m$. The proof follows parts of the proof of
\cite{Bedulev-Viehweg00, 3.1} very closely, however both the situation
and the statement are different from theirs, so the actual proof is
included.

\proclaim{\num Lemma}
Assume that $f$ is non-isotrivial and that $\omega_{X_{\text{gen}}}$
is ample. 

Then $(\fxc)^{\tensor e(m)r(m)}\tensor \det(\fxc)^{-1}$ is
semi-positive for all $m>0$.
\endproclaim\label\lemmafour

\proof
By \cite{Viehweg95, 2.8} one may replace $C$ by a finite cover,
unramified along $\Delta$, and $X$ by the pull-back family in order to
assume that $\det \fxc= \Cal D^{e(m)}$ for some invertible sheaf $\Cal
D$.

Let $r=r(m)$ and $\pi:Z\to X^r=X\times_C X\times_C\dots\times_C X$ a
resolution of singularities. Further let $\rho=f^r\circ \pi$.  Then
$\Cal M=\pi^*\omega_{X^r/C}=\pi^*\left(\otimes \text{pr}_i^*
\omega_{X/C} \right)$ is big by \cite{Viehweg83, Theorem II}. 

$f^r$ is a Gorenstein morphism and the general fibre has rational
singularities, so there are natural injective maps:
\vskip-.75cm
$$
\align
\rho_*\left(\Cal M^{m-1}\tensor\omega_{Z/C})\right) 
&\hookrightarrow
f^r_*\omega_{X^r/C}^m
\mytag\egystar
\\
f^r_*\omega_{X^r/C}^m
&\hookrightarrow
\rho_*\Cal M^m\mytag\ketstar
\\
\Cal D^{e(m)}
&\hookrightarrow (\fxc)^{\tensor r} \simeq
f^r_*\omega_{X^r/C}^m,\mytag\haromstar 
\endalign 
$$
where \quote\egystar\ and \quote\ketstar\ are isomorphisms near the
 generic point of $C$.

The composition of \quote\ketstar\ and \quote\haromstar\ gives a
section $\sigma\in\coh 0. Z. \Cal M^m\tensor \rho^*\Cal D^{-e(m)}.$. Let
$A=(\sigma=0)$. Since $\pi$ was an arbirary resolution of
singularities one may replace it by further blow-ups, so in particular
one may assume that $A$ is a normal crossing divisor.

Let $\Cal J\subseteq\ring Z.$ be the ideal sheaf defined as
$$
\im [ \rho^*f^r_*\omega_{X^r/C}^m \to \Cal M^m]=\Cal M^m\tensor \Cal J.
$$
Note that the support of $\ring X./\Cal J$ is contained in finitely
many fibers. By blowing up $\Cal J$ one can assume that it is a line
bundle and it is trivial near the general fibre of $\rho$. By
\cite{Kawamata82} $\fxc\simeq \gyc$ is semi-positive, hence so is $\Cal
M^m\tensor
\Cal J$, i.e., it is a nef line bundle.

Let $\Cal K=\Cal M^{m-1}\tensor \rho^*\Cal D^{-1}$. Then
$$
 \Cal K^{e(m)m}(-mA)=\Cal M^{e(m)m(m-1)-m^2}\supseteq (\Cal M^m\tensor
\Cal J)^{e(m)(m-1)-m},
$$
where the inclusion is an equality near the general fiber of
$\rho$. Hence $\Cal K^{e(m)}(-A)$ is nef near the general fiber of
$\rho$ and then $\rho_*\left(\Cal
K\tensor\omega_{Z/C}\left(-\left[\frac{A}{e(m)}\right]\right)\right)$
is semi-positive by \cite{Esnault-Viehweg90, 1.7}.  By
\cite{Viehweg95, 5.14, 5.21}
$$
 \Cal F:= \rho_*\left(\Cal
 K\tensor\omega_{Z/C}\left(-\left[\frac{A}{e(m)}\right]\right)\right)
 \hookrightarrow \rho_*(\Cal K\tensor\omega_{Z/C})
$$ 
is an isomorphism near the generic point of $C$.
On the other hand by \quote\egystar,
$$
\rho_*(\Cal K\tensor\omega_{Z/C}))\simeq \rho_*(\Cal M^{m-1}\tensor
\omega_{Z/C})\tensor\Cal D^{-1}
\hookrightarrow 
(\fxc)^{\tensor r}\tensor\Cal D^{-1},
$$
is also an isomorphism near the generic point of $C$.

Thus $\Cal F\subseteq (\fxc)^{\tensor r}\tensor\Cal D^{-1}$ which is
an equality on an open dense subset of $C$. $\Cal F$ is semi-positive
and then so is $(\fxc)^{\tensor r}\tensor\Cal D^{-1}$. Taking the
$e(m)$-th power gives the statement.
\endproof

\proclaim{\num Corollary}
Assume that $h$ is not birationally isotrivial and that $\omega_\ygen$
is nef and big. Then $(\gyc)^{\tensor e(m)r(m)}\tensor
\det(\gyc)^{-1}$ is semi-positive for all $m>0$.
\endproclaim\label\lemmafourb

\proclaim{\num Lemma}
Let $\Cal M$ be a line bundle on $X$ and $\Cal N$ a line bundle on
$C$. Assume that $f_*\Cal M\tensor \Cal N$ is ample,
$\Cal M_t=\Cal M\resto{X_t}$ is generated by global sections for $t\in
C\setminus\Delta$, and $h^0(\Cal M)$ is constant. Then

\myroster
\myitem
$\Cal M\tensor f^*\Cal N$ is semi-ample with respect to $X\setminus
D$.

\myitem
If $\Cal M_t$ is ample for $t\in C\setminus\Delta$, then $\Cal
M\tensor f^*\Cal N$ is ample with respect to $X\setminus D$.

\endmyroster

\endproclaim\label\lemmafive

\proof
Let $t,s\in C\setminus\Delta$. For $l\gg 0$, $\coh 1.C. \sym^l(f_*\Cal
M)\tensor \Cal N^l\tensor{\ring C.(-t-s)}.=0$. Hence the map
$
 \coh 0.C. \sym^l(f_*\Cal M)\tensor \Cal N^l. 
\to
 \sym^l(f_*\Cal M)\tensor (k(t)\oplus k(s))
$
is surjective. 

Since $\Cal M_t$ is generated by global sections for $t\in
C\setminus\Delta$, $f^*\sym^l(f_*\Cal M)\to \Cal M^l$ on $X\setminus
D$. In particular $ f^*\sym^l(f_*\Cal M)\tensor ({\ring X_t.}\oplus
{\ring X_s.})\to \Cal M_t^l\oplus \Cal M_s^l$ is surjective.

Now one has the following commutative diagram:
\diagram[size=.6cm]
 \coh 0.C. \sym^l(f_*\Cal M)\tensor \Cal N^l. \tensor {\ring X.}  &
\rTo^{\quad \alpha\quad } & f^*\sym^l(f_*\Cal M)\tensor ({\ring
X_t.}\oplus {\ring X_s.}) \\
\dTo & & \\
 \coh 0.C. f_*\Cal M^l\tensor \Cal N^l. \tensor {\ring X.} && \dTo_\beta \\
\dTo^\simeq && \\
 \coh 0.X. \Cal M^l\tensor f^*\Cal N^l. \tensor {\ring X.} &
 \rTo_\gamma &
\Cal M_t^l\oplus \Cal M_s^l
\enddiagram
Since $\alpha$ and $\beta$ are surjective, so is $\gamma$. This shows
both statements.
\endproof

\demo{\bf \num Definition}
Let $m_0(k)$ be the smallest positive constant \st\ for all projective
varieties, $F$, of dimension $k$ with at most rational singularities
and $\omega_F$ ample, $\omega_F^m$ is generated by global sections if
$m\geq m_0(k)$.
\enddemo\label\ggthreshhold

\demo{\remnum Remark}
Fujita's conjecture predicts that $m_0(k)\leq k+2$. 
\enddemo\label\fujita

\proclaim{\num Theorem}
Let $C$ be a smooth projective curve of genus $g$, $\Delta\subseteq C$
a finite set of points. Then there exists a divisor $\bar\Delta\subset
C$ of degree $2g+\delta+1$ \st\ for all non-isotrivial morphisms,
$f:X\to C$, satisfying the assumptions made in
\quote\assumptions\ and such  that for $t\in C\setminus \Delta$,
$\omega_{X_t}$ is ample, and 
for all $m\geq m_0(\dim X-1)$, 
$$
{\deg (f_*\omega_{X/C}^m)=\deg (\gyc)}\leq 
\deg(\omega_C(\bar\Delta)^{\dim X} 
\tensor\omega_C^{-1})^{me(m)r(m)}.  
$$
\endproclaim\label\bounded

\proof
$X$ has rational singularities, so $\phi_*\omega_{Y/C}\simeq
\omega_{X/C}$, and then $\gyc\simeq f_*\omega^m_{X/C}$.

{\bf\subnum} Let $m \geq m_0(\dim X-1)$ and $\Cal N$ a line bundle on
$C$
\st\
$$
\deg \Cal N^{-me(m)r(m)}<\deg (\fxc).
$$
Then by \quote\lemmafour\ $
f_*\omega^m_{X/C}\tensor \Cal N^{m}$ is ample on $C$. 
$h^0(X_t, \omega_{X_t}^m)$ is constant for $t\in C\setminus\Delta$ by
Kawamata-Viehweg vanishing, so $\omega_{X/C}\tensor f^* \Cal N$
is ample with respect to $X\setminus D$ by \quote\lemmafive.

{\bf\subnum} Choose an $l>0$ \st\ $\omega_{X/C}^l\tensor f^* \Cal
 N^{l}$ is generated by global section on $X\setminus D$. By blowing
 up the base locus of $\omega_{X/C}^l\tensor f^* \Cal N^{l}$
 (contained in $D$), one may assume that \te s an effective Cartier
 divisor $\Gamma$, supported on $\supp D$ \st\ $\omega_{X/C}^l\tensor
 f^* \Cal N^{l}(-\Gamma)$ is generated by global sections on the
 entire $X$.

{\bf\subnum} Let $P\in C\setminus \Delta$. We may assume that $l\geq
2g+\delta$. The linear system $|(2g+\delta)P-\Delta|$ is base point
free, so one can find a reduced effective divisor, $\Delta'\in
|(2g+\delta)P-\Delta|$
\st\ $\Delta\cap\Delta'=\emptyset$.  Let $D=\sum d_i D_i$ and
$l'=\underset{i}\to{\opname{lcm}}(d_i)\cdot l$.  Let
$$
\Delta''=\Delta + \Delta' + (l'-(2g+\delta))P \in |l'P|.
$$

Let $\sigma:\tilde C\to C$ be the finite cover obtained by taking the
$l'$-th root of $\Delta''$. Take the fiber product of $\sigma$ with
$f$ and $h$. Let $\tilde X$ be the normalization of $\tilde C\times_C
X$, $\tilde f:\tilde X\to \tilde C$, and $\tilde Y\to \tilde C\times_C
Y$ a resolution of singularities \st\ it is an isomorphism over
$\tilde C\setminus\sigma^*\Delta''$.

Note that $\Delta+\Delta'$ is a non-empty reduced divisor, so both
$\tilde C$ and $\tilde X$ are irreducible. Let
$\bar\Delta=(\Delta'')_{\text{red}}$,
$\tilde\Delta=(\sigma^*\bar\Delta)_{\text{red}}$, $\tilde D=\tilde
f^*\tilde \Delta$ and $\tilde D_j=(\pi^*D_j)_{\text{red}}$. Then
$\tilde Y\setminus \tilde\phi^*\tilde D\to \tilde
C\setminus\tilde\Delta$ is smooth, $\omega_{\tilde C}(\tilde
\Delta)\simeq\sigma^*\omega_C(\bar\Delta)$ and
$\bar\delta=\#\bar\Delta=2g+\delta+1$.  By \cite{Kov\'acs96; (2.4),
(2.17)} $\delta>0$ if $g=0$, so $2g-2+\bar\delta\geq 0$, i.e.,
$\omega_{\tilde C}(\tilde\Delta)\simeq\sigma^*\omega_C(\bar\Delta)$ is
nef. Furthermore, over the smooth locus of $X$, $
\pi^*D_j = \frac{l'}{\opname{gcd}(l, d_j)} \tilde D_j,
$ so by the definition of $l'$, the coefficient of $\tilde D_j$ is
divisible by $l$, hence \te s a divisor $\tilde \Gamma$ on $\tilde X$,
supported on $\supp \tilde D$, \st\ $\pi^*\Gamma= l\,\tilde
\Gamma$. As before, by blowing up the ideal sheaf of $\tilde\Gamma$
one may assume that it is a Cartier divisor. Then
$$
\pi^*(\omega_{X/C}^l\tensor f^* \Cal N^{l}(-\Gamma))=
\pi^*(\omega_{X/C}\tensor f^* \Cal N)^l(-l\,\tilde\Gamma)=
(\pi^*(\omega_{X/C}\tensor f^* \Cal N)(-\tilde\Gamma))^l,
$$
so $\pi^*(\omega_{X/C}\tensor f^* \Cal N)(-\tilde\Gamma)$ is
semi-ample on $\tilde X$ and ample with respect to $\tilde
X\setminus\tilde D$. Finally let
$\tilde\omega=\pi^*\omega_{X/C}(-\tilde\Gamma)$ and $\tilde\Cal
N=\sigma^*\Cal N$. Using this notation $\tilde\omega\tensor\tilde
f^*\tilde \Cal N$ is semi-ample on $\tilde X$ and ample with respect
to $\tilde X\setminus\tilde D$.

{\bf\subnum} Let $\Cal K=\omega_C(\bar \Delta)^n$ where $n=\dim X$ and
let $\tilde\Cal K=\sigma^*\Cal K$. Then by construction
$$
\tilde\omega\tensor\tilde f^*(\tilde\Cal N\tensor\tilde\Cal
K)\subseteq\pi^*(\omega_{X/C}\tensor f^*(\Cal N\tensor\Cal K)).
$$
Let $\Cal L= \omega_{X/C}\tensor f^*(\Cal N\tensor\Cal K) \tensor
 f^*\omega_{C}(\bar\Delta)^{-1}$ and $\tilde\Cal
 L=\tilde\omega\tensor\tilde f^*\tilde\Cal N\tensor\tilde
 f^*\omega_{\tilde C}(\tilde \Delta)^{n-1}\subseteq \pi^*\Cal
 L$. Since $\omega_{\tilde C}(\tilde \Delta)$ is nef, $\tilde\Cal L$
 and $\tilde\Cal L\tensor
\tilde f^*\omega_{\tilde C}(\tilde \Delta)^{-(n-1)}$ are semi-ample on
$\tilde X$ and ample with respect to $\tilde X\setminus\tilde D$.
Hence {$\coh n.X.\omega_{X/C}\tensor f^*(\Cal N\tensor\Cal K).=
0$} by \quote\mainket.

Finally take $\Cal N=\omega_C(\bar\Delta)^{-n}\tensor\omega_C$. Then
$\Cal N\tensor \Cal K\simeq \omega_C$,and $\omega_{X/C}\tensor
f^*(\Cal N\tensor\Cal K)\simeq \omega_X$. Now $\coh n.X.\omega_X.\neq
0$, so $\deg (\fxc)\leq \deg\Cal N^{-me(m)r(m)}=
(\omega_C(\bar\Delta)^n\tensor\omega_C^{-1})^{me(m)r(m)}.$
\endproof

\proclaim{\num Corollary}
Under the assumptions of \quote\bounded, $2g-2+\delta>0$. In
particular \quote\weakbound\ holds.
\endproclaim\label\css

\proof 
Assume the contrary, i.e., either $g=0$ and $\delta\leq 2$ or $g=1$
and $\delta=0$.  First of all we may assume that $h:Y\to C$ is
semi-stable and in both cases \te s a finite endomorphism, $\tau:C\to
C$, of degree $>1$ such that $\tau$ is smooth over $C\setminus \Delta$
and completely ramified over $\Delta$. Hence $f_\tau:X_\tau\to C_\tau$
again has the same properties as $f$. In particular $C_\tau\simeq C$
and $h_\tau$ is smooth over $C_\tau\setminus \Delta_\tau\simeq
C\setminus \Delta$.  However,
$\deg({h_\tau}_*\omega^m_{Y_\tau/C_\tau})
=\deg\tau\cdot\deg({h}_*\omega^m_{Y/C})$, contradicting the
boundedness of $\deg({h}_*\omega^m_{Y/C})$.
\endproof

\proclaim{\num Corollary}
Under the assumptions of \quote\bounded,
$$
\frac{\deg (f_*\omega_{X/C}^m)}{\rank (f_*\omega_{X/C}^m)}=
\frac{\deg (\gyc)}{\rank(\gyc)}
\leq 
4\cdot \dim X\cdot(2g-2+\delta)\cdot{m\cdot e(m)}\qquad 
\text{for }m\geq {\dim X\choose 2} +2.
$$
\endproclaim\label\fourbound

\proof
By \quote\css\ $2g-2+\delta>0$, so 
$$
 \deg(\omega_C(\bar\Delta)^{\dim X} \tensor\omega_C^{-1})=\dim
X(4g-1+\delta)-(2g-2)\leq 4\dim X(2g-2+\delta).
$$
Furthermore $m_0(\dim X-1)\leq {\dim X\choose 2} +2$ by
\cite{Kolll\'ar97, 5.8} (\cf{Angehrn-Siu95}).
\endproof

\demo{\remnum Remark}
Fujita's conjecture suggests that \quote\fourbound\ should hold for
$m\geq \dim X+1$.
\enddemo

\proclaim{\num Corollary}
Let $h:Y\to C$ be a non birationally isotrivial family \st\ $Y$ is a
smooth projective variety of dimension $n$, $C$ is a smooth projective
curve of genus $g$ and \te s a finite subset $\Delta\subset C$ \st\
$h$ is smooth over $C\setminus \Delta$. Assume that $\omega_{Y/C}$ is
$h$-nef and $h$-big. Then 
$$
\frac{\deg (\gyc)}{\rank(\gyc)} \leq 
4\cdot \dim X\cdot(2g-2+\delta)\cdot{m\cdot
e({\omega^m_{Y_{\text{gen}}}})}\qquad
\text{for }m\geq {\dim X\choose 2} +2.
$$
\endproclaim\label\nefandbig

\proof
By the relative base point free theorem \cite{KMM87} \te\ morphisms
$\phi: Y\to X$ and $f:X\to C$ that satisfy \quote\assumptions, so
the statement follows from \quote\fourbound.
\endproof

\newnum 
Let $\frak D_h$ denote the moduli functor of canonically polarized
normal projective varieties with rational Gorenstein singularities and
Hilbert polynomial $h(t)$. Let $\frak D_h^{(m)}$ denote the submoduli
functor of varieties with $\omega_F^{m}$ very ample
\cf{Viehweg95, 1.20, 8.18}.
 
$\frak D_h^{(m)}$ is locally closed by \cite{Kawamata99} and is
bounded by definition, so by \cite{Viehweg95, 8.20} \te s a coarse
quasi-projective moduli scheme $D_h^{(m)}$ for $\frak D_h^{(m)}$.
Furthermore by \cite{Koll\'ar90} \te s an integer $p>0$ and a very
ample line bundle $\lambda$ on $D_h^{(m)}$ \st\ for any $f:T\to S\in
\frak D_h^{(m)}(S)$, $\lambda$ induces ${\det(f_*\omega_{T/S}^{m})}^p$
on $S$. Note also that by \cite{Viehweg95, 5.17} $e(\omega_F^m)$ is
bounded on $D_h^{(m)}$.

Let $\bar D_h^{(m)}$ be the projective closure of $D_h^{(m)}$
corresponding to the embedding given by $\lambda$, and $\bold
H=\hom((C, C\setminus \Delta),(\bar D_h^{(m)}, D_h^{(m)}))$ the scheme
parametrizing morphisms $\Psi: C\to \bar D_h^{(m)}$ \st\
$\Psi(C\setminus \Delta)\subset D_h^{(m)}$.
\endnewnum


\proclaim{\num Theorem}
There exists a subscheme of finite type $T\subset \bold H$ that
contains all points $[\Psi: C\to \bar D_h^{(m)}]\in \bold H$ induced
by morphisms $f:X\to C\in \frak D_h^{(m)}(C)$ \st\ $f\resto{X\setminus
f^*\Delta}$ admits a simultaneous resolution. $\square$
\endproclaim


\def\myquad{\ }


\enddocument